\documentclass[a4paper,12pt]{article}

\usepackage[latin1]{inputenc} 
\usepackage[T1]{fontenc}     

\usepackage{geometry}         
\usepackage[francais, english]{babel}  

\usepackage{amsmath}
\usepackage{amsfonts}
\usepackage{amssymb}
\usepackage{fancyhdr}
\usepackage{url}
\usepackage{mathrsfs}
\usepackage{euscript}

\usepackage[all]{xy}

\newtheorem{df}{Definition}[section]

\newtheorem{prop}[df]{Proposition}
\newtheorem{thm}[df]{Theorem}
\newtheorem{thrm}{Theorem}
\newtheorem{lem}[df]{Lemma}
\newtheorem{crl}[df]{Corollary}
\newcommand{\prf}{\noindent \textit{Proof}}
\newtheorem{rmk}[df]{Remark}
\newtheorem{expl}[df]{Example}
\newcommand{\dbar}{\overline{\partial}}
\newcommand{\ddbar}{\partial\overline{\partial}}

\newcommand{\R}{\mathbb{R}}

\newcommand{\C}{\mathbb{C}}
\newcommand{\Z}{\mathbb{Z}}
\newcommand{\N}{\mathbb{N}}

\newcommand{\ric}{\operatorname{Ric}}
\newcommand{\tr}{\operatorname{tr}}

\newcommand{\vol}{\operatorname{vol}}

\newcommand{\vareps}{\varepsilon}

\newcommand{\inj}{\operatorname{inj}}

\newcommand{\id}{\operatorname{id}}

\newcommand{\riem}{\operatorname{Rm}}
\newcommand{\e}{\mathbf{e}}
\newcommand{\f}{\mathbf{f}}
\newcommand{\F}{\mathbf{F}}

\newcommand{\diff}{\mathbf{Diff}}
\newcommand{\met}{\mathbf{Met}}
\newcommand{\cqfd}{ \hfill $\square$ }

\title{\textbf{From ALE to ALF gravitational instantons}}  
\author{\textsc{Hugues AUVRAY}}
\date{}

\begin{document}

\makeatletter
\renewcommand%
   {\section}%
   {%
   \@startsection{section}%
      {1}%
      {0mm}%
      {\baselineskip}%
      {0.5\baselineskip}%
      {\sc\large\centering}%
   }%
\makeatother

\makeatletter
\renewcommand%
   {\subsubsection}%
   {%
   \@startsection{subsubsection}%
      {1}%
      {0mm}%
      {1.25\baselineskip}%
      {0.25\baselineskip}%
      {\sf\normalsize}%
   }%
\makeatother

\makeatletter
\renewcommand%
   {\paragraph}%
   {%
   \@startsection{paragraph}%
      {4}%
      {0mm}%
      {0mm}%
      {-5pt}%
      {\it\normalsize}%
   }%
\makeatother

\maketitle

\renewcommand{\abstractname}{Abstract}
 \begin{abstract}
  We give an original analytic construction of hyperkähler ALF metrics on some ALE spaces of dihedral type, 
namely the spaces corresponding to minimal resolutions of $\C^2/\mathcal{D}_k$, with $\mathcal{D}_k$ the binary dihedral group of order $4k$ for an integer $k\geq2$.  
 \end{abstract}

 \begin{center}
   \rule{3cm}{0.25pt}
 \end{center}

~

\section*{Introduction}

The aim of this paper is to construct ALF gravitational instantons from the data of ALE ones, on minimal resolutions of $\C^2/\Gamma$, for $\Gamma$ some binary dihedral group. 
Let us start with a brief review on ALE gravitational instantons. 

~

\paragraph{Basic material.} 
Consider $\R^4$ identified to $\C^2$ by the standard canonical relation between coordinates, i.e.
 \begin{equation*}
  z_1=x_1+ i x_2, \quad z_2=x_3+ix_4.
 \end{equation*}
One shall denote the associated complex structure by $I_1$; we also denote respectively by $I_2$ and $I_3$ the complex structures associated to the complex coordinates $x_1+ i x_3$ and $x_4+ i x_2$ on the one hand, and $x_1+ i x_4$ and $x_2+ i x_3$ on the other hand. 
This way, one has the quaternionic relations $I_1I_2=-I_2I_1=I_3$, etc.

Let us take $I_1$ as a favorite complex structure, meaning that without further precisions, we refer implicitly 
to this complex structure. 
We define a family of finite subgroups of $SU(2)$ acting on $\R^4$, the \textit{binary dihedral} groups: 
 \begin{df}[Binary dihedral groups.] \label{df_dk}
  For $k\geq2$, let $\zeta_k$ act on $\R^4\simeq \C^2$ \textit{via} $\zeta_k(z_1,z_2)=\big(e^{\tfrac{i\pi}{k}}z_1,e^{-\tfrac{i\pi}{k}}z_2\big)$, and $\tau$ \textit{via} $\tau(z_1,z_2)=(z_2, -z_1)$. 
We call the group $\langle \zeta_k, \tau\rangle$ the dihedral group of order $4k$, and denote it by $\mathcal{D}_{k}$.
 \end{df}

 \begin{rmk}
  We are focusing in what follows on the \textit{finite groups} of this definition. 
  In order to avoid confusion though, let us note that to these finite groups correspond some Lie groups and root systems.  
  Notice the following shift: to the finite group $\mathcal{D}_k$ corresponds the Lie groups "of type $D_{k+2}$". 
  This being precised, when using the notation $\mathcal{D}_k$, we will always refer to the finite group; one can thus forget about the Lie groups alluded here.
 \end{rmk}

If $\Gamma$ is one of the groups above, then the standard euclidean metric $\e$ on $\C^2$ is clearly invariant under the action of $\Gamma$ which acts freely outside of $0$, hence descends to $(\R^4\backslash\{0\})/\Gamma$. 
For the same reason, $I_1$ (but also $I_2$ and $I_3$) and $r$ make sense on $(\R^4\backslash\{0\})/\Gamma$; throughout the paper, we denote the projection $\R^4\to\R^4/\Gamma$ by $\pi$. 

Let us fix $k\geq 2$. 
Notice that the three polynomials 
 \begin{equation*}
  u=\frac{1}{2}\big(z_1^{2k+1}z_2-z_2^{2k+1}z_1\big),\quad v=\frac{i}{2}\big(z_1^{2k}+z_2^{2k}\big),\quad w=z_1^2 z_2^2
 \end{equation*}
are invariant under the action of $\mathcal{D}_k$, and verify the relation: $u^2+v^2w+w^{k+1}=0$. 
Actually, one may also check that $\C[z_1,z_2]^{\mathcal{D}_k}=\C[u,v,w]$, and moreover that 
 \begin{equation*}
  X_{\mathcal{D}_k}:=\{(u,v,w)\in\C^3|\, u^2+v^2w+w^{k+1}=0\} \simeq  \C^2/\mathcal{D}_k
 \end{equation*}
as complex orbifolds. 

~

Now, ALE gravitational instantons can be defined as follows: 
 \begin{df}
  We say that $(X,g)$ is an Asymptotically Locally Euclidean gravitational instanton, or ALE gravitational instanton for short, if 
  it is a hyperkähler complete manifold of real dimension 4, of maximal volume growth, and with square integrable Riemannian curvature.
 \end{df} 
Notice that the Ricci-flatness of such spaces follows at once from the hyperkähler assumption. 
From the works \cite{bkn} and \cite{kro2}, those spaces are completely classified 
-- let us mention quickly a recent addition to this classification: 
if one drops the hyperkähler assumption (keeping the Ricci-flatness, the Kähler assumption and the ALE geometry), 
one has to add quotients of certain spaces from Kronheimer's list \cite{suv}. 
We will not give here this classification, but instead underline some of its consequences. 

To an ALE gravitational instanton $(X,g)$ is attached a finite subgroup $\Gamma$ of $SU(2)$; then : 
 \begin{itemize}
  \item[$\bullet$] $X\backslash K$ is diffeomorphic to $(\R^4\backslash B)/\Gamma$, where $K$ is some compact subset of $X$, $B$ is a ball of $\R^4$ centered at the origin; 
  \item[$\bullet$] let $\underline{g}$ be the pull-back to the metric to $(\R^4\backslash B)$; then $\nabla^{\ell}_{\e}(\e-\underline{g})=O(r^{-4-\ell})$ for all $\ell\geq 0$, provided that the diffeomorphism above is well-chosen ;
  \item[$\bullet$] choose three complex structures $I^X_{j}$ on $X$ such that the data $(g,I^X_{j})$ is a hyperkähler structure on $X$; 
                   then up to the action of $SO(3)$ on the $I^X_{j}$, one has $\nabla^{\ell}_{\e}(I_j-\underline{I^X_{j}})=O(r^{-4-\ell})$ for all $\ell\geq 0$. 
                   Here the underlining indicates the pulling-back to $\C^2$, for the \textit{same} diffeomorphism as that of the previous point.
 \end{itemize}
 \begin{expl}[Eguchi-Hanson metric] \label{expl_eh}
  The Eguchi-Hanson metric on the minimal resolution (for the complex structure $I_1$) $T^*P^1_{\C}$ of $\C^2/\Z_2$ is a standard example. 
  If one identifies $(\C^2\backslash\{0\})/\Z_2$ and $T^*P^1_{\C}\backslash P^1_{\C}$, then the (pulled-back) metric can be written as $\e+g_1+\vareps$, with :
   \begin{align*}
    g_1=&\frac{(|z_1|^2-|z_2|^2)}{r^6}(dx_1^2+dx_2^2-dx_3^2-dx_4^2)+\frac{2(x_1x_3+x_2x_4)}{r^6}(dx_1\cdot dx_3+dx_2\cdot dx_4)\\
        &+\frac{2(x_1x_4-x_2x_3)}{r^6}(dx_1\cdot dx_4-dx_2\cdot dx_3),
   \end{align*}
  where $\alpha\cdot\beta=\alpha\otimes\beta+\beta\otimes\alpha$, and
   \begin{equation*}
    \big|(\nabla^{\e})^{\ell}\vareps\big|_{\e} = O(r^{-6-{\ell}})\quad \text{for all }\ell\geq0.
   \end{equation*}
  Here our identification plays the role of the diffeomorphism; the diffeomorphism is thus holomorphic for $I_1$ and $I_1^X$ in this situation. 
  However, one can only choose $I^X_2$ and $I^X_3$ so that the diffeomorphism is asymptotically holomorphic for the pairs $(I_2,I_2^X)$ and $(I_3,I_3^X)$, in the sense given above. 
  Notice also that $\tr^{\e}(g_1)=0$, so that actually $\det^{\e}(\e+g_1+\vareps)=1+O(r^{-6})$.
 \end{expl}

\paragraph{Results.} 
We are interested more specifically in this paper in ALE gravitational instantons $X$ of dihedral type, i.e. those whose attached group $\Gamma$ is one of the $\mathcal{D}_k$ of Definition \ref{df_dk}. 
In a nutshell, our work consists in constructing analytically metrics with different asymptotic geometry on (some of) those spaces, namely \textit{ALF metrics}, modeled on the Taub-NUT metric, described in next part -- let us only mention here that this metric is \textit{Ricci-flat} and has \textit{cubic volume growth}--, 
pushed-forward to $\C^2/\Gamma$; 
the acronym ALF stands here for \textit{Asymptotically Locally Flat}. 
The main result of this paper is the following, for $(X,g)$ an ALE gravitational instanton on a minimal resolution of $\C^2/\Gamma$:
 \begin{thrm}  \label{thrm_gRF}
  There exists a one-parameter family of hyperkähler metrics $g_{RF,m}$ on $X$, with volume form $\Omega_g$, 
  which are Kähler for $I_1^X$, indexed by a "mass" parameter $m\in(0,+\infty)$. 
  When $m$ is fixed, $(g_{RF,m}-\tilde{\f}_m)$ and $\nabla^{\tilde{\f}_m}(g_{RF,m}-\tilde{\f}_m)$ are $O(\tilde{R}^{-2-\delta})$ 
  for any $\delta\in (0,1)$. 
  If one denotes by $\omega_{RF,m}$ the Kähler form $g_{RF,m}(I_1^X\cdot,\cdot)$, 
  then $\omega_{RF,m}$ and $\omega_g$ are in the same cohomology class (Theorem \ref{thm_gRF}). 
 \end{thrm}

In this statement, $\tilde{\f}_m$ is a metric on $X$ equal to the push-forward of the Taub-NUT metric of parameter $m$ near infinity, 
and $\tilde{R}$ plays the role of the distance to some point for $\tilde{\f}_m$. 

Theorem \ref{thrm_gRF} is proved by gluing a potential of the Taub-NUT metric outside a compact subset of $X$, so as to obtain a Kähler metric (for the favorite complex structure), which has a small Ricci curvature near infinity. 
One then has to correct the Ricci-form in order to get a Ricci-flat Kähler metric on the whole $X$; 
to this aim we use a Calabi-Yau type theorem, which we state:
 \begin{thrm} \label{thrm_CY_ALF}
  Let $(Y,g_Y,J^Y,\omega_Y)$ be an ALF Kähler 4-manifold of dihedral type. 
  Given a $C^{\infty}_{loc}$ function $f$ such that there exists some positive $\beta\in(0,1)$ such that $(\nabla^{g_Y})^{\ell} f=O(\rho^{-2-\ell-\beta})$ for all $\ell\geq 0$, 
  there exists a unique $\varphi\in C^{\infty}_{loc}$ such that $(\nabla^{g_Y})^{\ell} \varphi=O(\rho^{-\ell-\beta})$ for all $\ell\geq 0$ and 
  \begin{equation*} 
   \big(\omega_Y+i\ddbar\varphi\big)^2=e^f\omega_Y^2
  \end{equation*}
 (Theorem \ref{thm_CY_ALF}).
 \end{thrm}
In this statement, $\rho$ is a distance function for $g_Y$. 
As is said above, an "ALF Kähler 4-manifold of dihedral type" is modeled on the Taub-NUT metric pushed-forward to $\C^2/\Gamma$; 
we are still vague about the meaning of this assertion, and will be clearer when using Theorem \ref{thrm_CY_ALF}. 

~

\paragraph{Organization of the article.} 
In part \ref{part_TN}, we give an original description of a 1-parameter family of Taub-NUT metrics on $\C^2$ endowed with its usual complex structure $I_1$; this description is suggested by LeBrun in \cite{leb}. 
In particular, we give a rather explicit potential for any metric of this family in section \ref{prgrph_taubnut_C2}. 
We also compare it to the standard euclidean metric $\e$; we specify a mutual control between them, and express objects attached to euclidean geometry (such as the $dx_j$) with the help of 1-forms and vector fields better adapted to the Taub-NUT metric (section \ref{sect_comp_ef}). 
Such a comparison is motivated by the computational considerations of the following parts. 

In part \ref{part_appl}, we apply the description of part \ref{part_TN} to the construction stated in Theorem \ref{thm_gRF}. 
A gluing of the appropriate potential is performed in section \ref{sect_gluing}. 
On the other hand, the correction of the Ricci form we alluded above needs several steps of analytic work, and is dealt with through sections \ref{sctn_RFnrinfty} to \ref{sctn_ccl_pgrm}. 

The conclusion of this program relies on Theorem \ref{thrm_CY_ALF}, which is proved in part \ref{part_prf_CY_ALF}, with a continuity method. 
The method is explained in section \ref{sect_cntnty_mthd}; the analysis it requires to work is  done in the following section, and the conclusion is given in the last section. 

~

Let us say now a few words about our construction. 
First, the objects we construct might not be new, see the constructions in \cite{ch-hi} and \cite{ch-ka}; 
though, as illustrated by Theorem \ref{thrm_gRF} and the idea of its proof we give, our analytic construction is quite explicit, and the models it relies on are quite close to the sought objects, which may ease their comprehension. 

We moreover restrict to minimal resolutions; we do not explore the cases of other types of finite subgroups of $SU(2)$, namely the cyclic groups, the tetrahedral group, the octahedral group and the dodecahedral group, and neither the general cases of dihedral ALE instantons.
This because:
 \begin{itemize}
  \item[$\bullet$] our following construction will not be possible for tetrahedral, octahedral and dodecahedral groups (we will use Kähler potentials invariant under the actions of cyclic or dihedral groups, but not under those of tetrahedral, octahedral and dodecahedral groups);
  \item[$\bullet$] ALF spaces of cyclic types are very well understood: their classification has been completely established by Minerbe \cite{min2}, and these spaces are explicitly described by the Gibbons-Hawking's ansatz; 
  \item[$\bullet$] in the general $\mathcal{D}_k$ cases, one has to deal with \textit{perturbed} complex structures, and the Kähler potentials we use behave badly with the error terms coming from this approximation. 
                   Though we do not circumvent this issue so far, we can at least focus on the case where one of the complex structures of the instanton agrees with one of those on $\C^2$. 
                   This corresponds to the situation where $X$ is a \textit{minimal resolution} of one of the orbifolds $X_{\mathcal{D}_k}$, the general case corresponding to $X$ being a \textit{deformation} of a $X_{\mathcal{D}_k}$, and 
                   we hope to deal with this general case in a future paper. 
                   The main difficulty here is to find a better potential for the rough metric, after what our method would apply globally unchanged; see section \ref{rmk_alpha0} for further details on this point.  
 \end{itemize}
Besides, a reading of part \ref{part_appl} may allow one to consider our gluing could have been better, e.g. using more closely algebraic properties of resolutions to get rid of some annoying error terms appearing in the process. 
Our choice of being non accurate at that point is nonetheless guided by the thought that the analytic work we lead in part \ref{part_appl} would be rather necessary in a similar construction on deformations of $\C/\Gamma$.

Let us underline now the following about Theorem \ref{thrm_gRF}. 
First, the expected behaviour of the family $(g_{RF,m})$ when $m$ goes to 0 is a convergence in 
$C^{\infty}_{loc}$ topology to the initial ALE metric $g$. This is quite obvious on the model $\C^2$, 
but because of our corrections, especially when using Theorem \ref{thrm_CY_ALF}, it is not clear yet on the considered spaces.
Notice  also that we do not restrict to presumably "small" deformations of the initial ALE metric into ALF metrics. 
This would correspond to work only with $m\in(0,m_0]$ for some finite $m_0$; instead we handle our construction \textit{for all} $m\in(0,+\infty)$.

To conclude, we comment briefly Theorem \ref{thrm_CY_ALF}. 
This result comes within the general scope of generalizing the celebrated Calabi-Yau theorem (see \cite{yau}) 
to non-compact manifolds, initiated by Tian and Yau \cite{tian-yau1,tian-yau2}. 
A statement similar to ours can be found in \cite[Prop. 4.1]{hein} (take the parameter $\beta$ of this statement equal to 3), in a more general and abstract framework. 
One interest of our result nonetheless, simply lies in the fact that besides we ask more precise asymptotics on our data than Hein does,
we get in compensation sharper asymptotics on our solution $\varphi$, which are useful in establishing Theorem \ref{thrm_gRF}.

~

\paragraph{Acknowledgements.} 
I am very grateful to my PhD adviser Olivier Biquard for introducing me to this subject of differential geometry, and for his precious help on this work. 
I also wish to thank Michael Singer for reading a previous version on this text, and Vincent Minerbe for indicating me C. LeBrun's paper \cite{leb}.

~

\section{Taub-NUT metric on $\C^2$}  \label{part_TN}

  The Taub-NUT metric on $\C^2$, introduced by Hawking \cite{haw}, is a general prototype for metrics of ALF gravitational instantons, which are roughly speaking 
complete hyperkähler 4-manifolds of cubic volume growth and whose Riemannian curvature has faster than quadratic decay; precisions are made below. 

To describe such a model, we solve an exercise in \cite{leb}, which is stated so as to provide a potential for the Taub-NUT metric for the standard complex structure $I_1$ on $\C^2$.

We then compare the Taub-NUT and the euclidean metrics. 

  \subsection{A potential for the Taub-NUT metric on $\C^2$}  \label{prgrph_taubnut_C2}

In \cite{leb} C. LeBrun suggests the following exercise: let $m$ be a positive parameter, and $u$ and $v$ implicitly defined on $\C^2$ by formulas:
 \begin{equation}  \label{eqn_formula_uv}
  \begin{aligned}
   |z_1|&=e^{m(u^2-v^2)}u, \\
   |z_2|&=e^{m(v^2-u^2)}v 
  \end{aligned}
 \end{equation}
(we do not make the dependence on $m$ appearent here, since for now we see this parameter as fixed; 
we shall only add $m$ as an index if we need to emphasize this dependence). 
One has: 
 \begin{prop}[LeBrun] \label{prop_leb}
  The metric $\f$ associated to the form 
   \begin{equation*}
    \omega_{\f}:=\frac{1}{4}dd^c \big(u^2+v^2+m(u^4+v^4)\big)
   \end{equation*}
  for the standard complex structure $I_1$ on $\C^2$ is the Taub-NUT metric.
 \end{prop}

Recall briefly the following description for the Taub-NUT metric: 
it takes the form $V(dy_1^2+dy_2^2+dy_3^2)+V^{-1}\eta^2$, where the $y_j$ are the coordinates of $\R^3$, 
$V=2m\big(1+\tfrac{1}{4mR}\big) $ where $R^2=y_1^2+y_2^2+y_3^2$, 
$\eta$ is a connection 1-form on a circle bundle over $\R^3\backslash\{0\}$ such that $d\eta=*_{\R^3}dV$; 
a complex structure is moreover given for example by $Vdy_1\mapsto \eta$, $dy_2\mapsto dy_3$, the other two being given by cyclic permutation of the $y_j$; see e.g. \cite{leb} for details.

\begin{rmk} \label{rmk_dk_invce}
 Notice that with the notations of Definition \ref{df_dk}, $u$ and $v$ are invariant under the action of $\zeta_k$, whereas that of $\tau$ interchanges them. 
As a consequence, $u^2+v^2+m(u^4+v^4)$ is invariant under $\mathcal{D}_k$, and is smooth on $(\C^2\backslash\{0\})/\mathcal{D}_k$. 
It is also invariant under the binary cyclic groups, but not anymore under the binary tetrahedral group 
-- generated by $\mathcal{D}_2$ and $\tfrac{1}{\sqrt{2}}\Big(\begin{smallmatrix} \vareps^7 & \vareps^7 \\ \vareps^5 & \vareps \end{smallmatrix}\Big)$, 
with $\vareps=e^{i\pi/4}$ --
the binary octahedral group -- generated by the binary tetrahedral group and $\Big(\begin{smallmatrix} \vareps & 0 \\ 0 & \vareps^7\end{smallmatrix}\Big)$ --
nor the binary icosahedral group -- generated by $ \Big(\begin{smallmatrix} -\eta^3 & 0 \\ 0 & -\eta^2\end{smallmatrix}\Big)$ 
and $ \tfrac{1}{\eta^2-\eta^3}\Big(\begin{smallmatrix} \eta+\eta^4 & 0 \\ 0 & -(\eta+\eta^4)\end{smallmatrix}\Big)$, with $\eta=e^{2i\pi/5}$.
\end{rmk}

As LeBrun leaves the proof to his reader, we shall give our own one here, relying on a sequence of basic lemmas. 
Before this, we shall mention that LeBrun's potential may be obtained by hyperkähler quotient considerations; 
we chose to give a less conceptual proof though since it exhibits several objects we shall need further.

 \begin{lem} \label{lem_f_ric_flat}
  The metric $\f$ is Ricci-flat; more precisely, $\omega_{\f}^2=2\Omega_{\e}$, where we recall that $\Omega_{\e}$ is the standard volume form $\tfrac{idz_1\wedge d\overline{z_1}\wedge idz_2\wedge d\overline{z_2}}{4}$.
 \end{lem}
\prf. We start by the computation of $\omega_{\f}$, which goes through that of $\tfrac{\partial u}{\partial z_j}$, $\tfrac{\partial v}{\partial z_j}$, $j=1,2$. 
One has: 
  \begin{equation}  \label{eq_deriv1}
   \begin{aligned}
    &\frac{\partial u}{\partial z_1}=\frac{1+2mv^2}{(2z_1)\big(1+2m(u^2+v^2)\big)}u,&   &\frac{\partial u}{\partial z_2}=\frac{muv^2}{z_2\big(1+2m(u^2+v^2)\big)},     \\
    &\frac{\partial v}{\partial z_1}=\frac{mu^2v}{z_1\big(1+2m(u^2+v^2)\big)},&         &\frac{\partial v}{\partial z_2}=\frac{1+2mu^2}{(2z_2)\big(1+2m(u^2+v^2)\big)}v.
   \end{aligned}
  \end{equation}
Indeed, differentiating the relation $|z_1|=e^{m(u^2-v^2)}u$ with $\tfrac{\partial }{\partial z_1}$, it follows that
 \begin{equation*}
  \frac{1}{2}\frac{|z_1|}{z_1}=\Big[m\Big(2u^2\frac{\partial u}{\partial z_1}-2uv\frac{\partial v}{\partial z_1}\Big)+\frac{\partial u}{\partial z_1}\Big]e^{m(u^2-v^2)},
 \end{equation*}
hence, writing $e^{m(u^2-v^2)}=\frac{|z_1|}{u}$, $u=2z_1\big[(1+2mu^2)\tfrac{\partial u}{\partial z_1}-2muv\tfrac{\partial v}{\partial z_1}\big]$. 
Similarly, applying $\tfrac{\partial }{\partial z_1}$ to the relation $|z_2|=e^{m(v^2-u^2)}v$, one gets $0=(1+2mv^2)\tfrac{\partial v}{\partial z_1}-2muv\tfrac{\partial u}{\partial z_1}$, that is $\tfrac{\partial v}{\partial z_1}=\tfrac{2muv}{1+2mv^2}\tfrac{\partial u}{\partial z_1}$. 
Substituting in the previous equality, one has
 \begin{equation*}
  \frac{\partial u}{\partial z_1}=\frac{1+2mv^2}{(2z_1)\big(1+2m(u^2+v^2)\big)}u,
 \end{equation*}
and then $\tfrac{\partial v}{\partial z_1}=\tfrac{mu^2v}{z_1(1+2m(u^2+v^2))}$. 
The other two equalities are obtained by symmetry.

Set now $\varphi=\tfrac{1}{4}\big(u^2+v^2+m(u^4+v^4)\big)$ (we will keep this notation throughout the paper). 
According to formulas \eqref{eq_deriv1}, 
 \begin{equation*}
  2\frac{\partial \varphi}{\partial z_1}=u(1+2mu^2)\frac{\partial u}{\partial z_1}+v(1+2mv^2)\frac{\partial v}{\partial z_1}=\frac{(1+2mv^2)u^2}{2z_1}
 \end{equation*}
and $2\tfrac{\partial \varphi}{\partial z_2}=\tfrac{(1+2mu^2)v^2}{2z_2}$, i.e. 
 \begin{equation*}
  \frac{\partial \varphi}{\partial \overline{z_1}}=\frac{(1+2mv^2)u^2}{4\overline{z_1}} \qquad \text{and} \qquad \frac{\partial \varphi}{\partial \overline{z_2}}=\frac{(1+2mu^2)v^2}{4\overline{z_2}}
 \end{equation*}
by conjugation. 
Apply again 
$\tfrac{\partial }{\partial z_1}$ and $\tfrac{\partial }{\partial z_2}$ to those equalities, as well as the relation $uv=|z_1z_2|$ and the formulas \eqref{eq_deriv1}, it follows, setting $R=\tfrac{1}{2}(u^2+v^2)$:
 \begin{align*}
  \omega_{\f} =&dd^c\varphi\\
           =& \Big(\frac{u^2(1+2mv^2)}{2|z_1|^2(1+4mR)}+m|z_2|^2\Big)idz_1\wedge d\overline{z_1}
             +m\overline{z_2}z_1\Big(1+\frac{1}{1+4mR}\Big)idz_2\wedge d\overline{z_1}\\
            & +m\overline{z_1}z_2\Big(1+\frac{1}{1+4mR}\Big)idz_1\wedge d\overline{z_2}       
              +\Big(\frac{v^2(1+2mu^2)}{2|z_2|^2(1+4mR)}+m|z_1|^2\Big)idz_2\wedge d\overline{z_2} .  
 \end{align*}
The conclusion follows from direct computation of $\omega_{\f}^2$, using $uv=|z_1z_2|$ once more. \cqfd

~

Let $S^1=\R/2\pi\Z$ act on $\C^2$ by $\alpha\cdot(z_1,z_2)=(e^{i\alpha}z_1,e^{-i\alpha}z_2)$. 
The infinitesimal action is given by the vector field 
 \begin{equation*}  \label{formule_xi}
  \xi:=i\Big(z_1\frac{\partial }{\partial z_1}+\overline{z_2}\frac{\partial }{\partial \overline{z_2}}-z_2\frac{\partial }{\partial z_2}-\overline{z_1}\frac{\partial }{\partial \overline{z_1}}\Big).
 \end{equation*}
By invariance of $u$ and $v$ under this circle action, it is clear that $\xi\cdot u=\xi\cdot v=\xi\cdot \varphi=0$, and similarly that $\mathcal{L}_{\xi}\omega_{\f}=0$. 
Let us furthermore consider the $(2,0)$-form $\Theta:=dz_1\wedge dz_2$, which is symplectic and holomorphic, so that $\Theta\wedge\overline{\Theta}=4\Omega_{\e}=2\omega_{\f}^2$. 
One moreover has $\mathcal{L}_{\xi}\Theta=0$.

We look for a complex hamiltonian $H$ for the $S^1$-action, relatively to the symplectic structure given by $\Theta$. 
Now one has $\iota_{\xi}\Theta=(z_1dz_2 + z_2 dz_1)=d(iz_1z_2)$, which prompts us to set $y_2:=\mathfrak{Im}(z_1z_2)$ and $y_3:=-\mathfrak{Re}(z_1z_2)$, i.e. $H=y_2+iy_3$.

In the same way, one has  $\mathcal{L}_{\xi}d^c\varphi=0$, and the right-hand-side member also equals, by Cartan's formula, $\iota_{\xi}dd^c\varphi+d(\iota_{\xi}d^c\varphi)$, hence: $\iota_{\xi}\omega_{\f}=-d\big(d^c\varphi(\xi)\big)$, which encourages us to set $y_1=d^c\varphi(\xi)$ this time. 
All computations done, we state:
 \begin{lem}
  One has $y_1=\tfrac{1}{2}(u^2-v^2)$, and thus $R$ indeed equals $(y_1^2+y_2^2+y_3^2)^{1/2}$.
 \end{lem}
\prf. According to the proof of Lemma \eqref{lem_f_ric_flat}, 
 \begin{equation*}
  d^c\varphi=i(1+2mv^2)u^2\Big(\frac{d\overline{z_1}}{2\overline{z_1}}-\frac{dz_1}{2z_1}\Big)+i(1+2mu^2)v^2\Big(\frac{d\overline{z_2}}{2\overline{z_2}}-\frac{dz_2}{2z_2}\Big),
 \end{equation*}
hence the result from formula \eqref{formule_xi}. 

To get the equality $y_1^2+y_2^2+y_3^2=\tfrac{1}{4}(u^2+v^2)^2$, it is enough to notice that $y_2^2+y_3^2=|y_2+iy_3|^2=|z_1z_2|^2=u^2v^2$. 
\cqfd

~

Finally, let us set $V=|\xi|^{-2}_{\f}$. 
 \begin{lem}
  One has $|\xi|^{2}_{\f}=\tfrac{2R}{1+4mR}$, and hence $V=2m\big(1+\tfrac{1}{4mR}\big)$.
 \end{lem}
\prf. As $I_1$ is the standard complex structure of $\C^2$, $I_1\xi=-z_1\tfrac{\partial }{\partial z_1}-\overline{z_1}\frac{\partial }{\partial \overline{z_1}}+\overline{z_2}\tfrac{\partial }{\partial \overline{z_2}}+z_2\frac{\partial }{\partial z_2}$. 
The easy but tedious calculation of $|\xi|_{\f}^2=\omega_{\f}(\xi,I_1\xi)$ then follows. 
However, noticing  that
 \begin{align*}
  &idz_1\wedge d\overline{z_1}(\xi,I_1\xi)=2|z_1|^2,&              &idz_1\wedge d\overline{z_2}(\xi,I_1\xi)=-2z_1\overline{z_2}, \\
  &idz_2\wedge d\overline{z_1}(\xi,I_1\xi)=-2z_2\overline{z_1},&   &idz_2\wedge d\overline{z_2}(\xi,I_1\xi)=2|z_2|^2;
 \end{align*}
makes easier the computation of $\sum_{j,k=1}^2 (\omega_{\f})_{j\bar{k}}idz_j\wedge d\overline{z_k}(\xi,I\xi)$. 
\cqfd

~

In order to get the Taub-NUT metric back under its classical form, we need finally a 1-form $\eta$, which would be a connection 1-form for the circle fibration $\varpi:\C^2\backslash\{0\}\mapsto \R^3\backslash\{0\}$, $(z_1,z_2)\mapsto (y_1,y_2,y_3)$. 
The natural candidate is given by 
 \begin{equation*}
  \eta:=VI_1dy_1.
 \end{equation*}

 \begin{lem}  \label{lem_eta}
  On $\C^2\backslash\{z_1z_2=0\}$, one has 
   \begin{equation*}
    \eta=\frac{i}{4R}\Big[u^2\Big(\frac{d\overline{z_1}}{\overline{z_1}}-\frac{dz_1}{z_1}\Big)-v^2\Big(\frac{d\overline{z_2}}{\overline{z_2}}-\frac{dz_2}{z_2}\Big)\Big], 
   \end{equation*}
  and $\eta(\xi)=1$ outside of 0.
 \end{lem}
\prf. By definition, $\eta=Vd^c y_1=\tfrac{1}{2}iV\big(2u(\dbar u-\partial u)-2v(\dbar v-\partial v)\big)$. 
We then apply formula \eqref{eq_deriv1} which gives the first order derivatives of $u$ and $v$, which can be rewritten (with the introduced notations) :
 \begin{align*}
   V\frac{\partial u}{\partial z_1}=\frac{1+2mv^2}{4z_1R}u,  \quad 
   V\frac{\partial u}{\partial z_2}=\frac{muv^2}{2z_2R},     \quad
   V\frac{\partial v}{\partial z_1}=\frac{mu^2v}{2z_1R},     \quad
   V\frac{\partial v}{\partial z_2}=\frac{1+2mu^2}{4z_2R}v,
 \end{align*}
hence (for instance), the component of $\eta$ in the direction of $dz_1$ is $-iu\frac{1+2mv^2}{4z_1R}u+iv\frac{mu^2v}{2z_1R}=-\tfrac{iu^2}{4z_1R}$, and similarly for the other three components. 
A straightforward computation suffices to see that $\eta(\xi)=1$. \cqfd

~

We shall now recover the Taub-NUT metric under a more familiar shape:
 \begin{lem}
  Away from 0, $\omega_{\f}= dx_1\wedge\eta+Vdx_2\wedge dx_3$, and thus $\f=V(dy_1^2+dy_2^2+dy_3^2)+V^{-1}\eta^2$.
 \end{lem}
\prf. This is rather clear that the family formed by $dy_2$ and $dy_3$ is linearly independent away from $0\in\C^2$, and that those forms vanish against $\xi$ because $y_2$ and $y_3$ are invariant under the action of $S^1$. 
We see that they vanish as well against $I_1\xi$, because for instance $I_1dy_2=dy_3$. 
Since $dy_1$ vanishes against $\xi$ ($y_1$ is $S^1$-invariant) but not against $I_1\xi$ since $I_1dy_1=V^{-1}\eta$, and since $\eta$ vanishes $I_1\xi$ but not against $\xi$, we deduce that $\{dy_1,dy_2,dy_3,\eta\}$ is free outside of $0$.

Consequently, on $\C^2\backslash\{0\}$, one has at any point the writing:
 \begin{align*}
  \omega_{\f}=\alpha dy_1\wedge\eta +\beta dy_2\wedge\eta  +\gamma dy_3\wedge\eta 
           +\delta dy_1\wedge dy_2+\vareps dy_1\wedge dy_3 +\zeta dy_2\wedge dy_3
 \end{align*}
for some $\alpha,\dots,\zeta$ (depending on the point). 
Now $\alpha  dy_1 +\beta  dy_2+\gamma dy_3=-\iota_{\xi}\omega_{\f}=dy_1$, thus $\alpha=1$ and $\beta=\gamma=0$, and since $\omega$ is of type $(1,1)$, one also has $\delta=\vareps=0$. 

It remains to determine the coefficient $\zeta$; we easily see that 
 \begin{equation*}
  (dy_1\wedge\eta)\Big(\frac{\partial}{\partial z_1},\frac{\partial}{\partial \overline{z_1}}\Big) = 2iV\Big|\frac{\partial y_1}{\partial z_1}\Big|^2 
                                                                                                   = i\frac{Vu^4}{2|z_1|^2(1+4mR)^2}
                                                                                                   = i\frac{u^4}{4R|z_1|^2(1+4mR)},
 \end{equation*}
whereas 
 \begin{equation*}
  \omega_{\f}\Big(\frac{\partial}{\partial z_1},\frac{\partial}{\partial \overline{z_1}}\Big)=i\frac{u^2(1+2mv^2)}{2|z_1|^2(1+4mR)}+im|z_2|^2
 \end{equation*}
hence
$\big(\omega_{\f}- dy_1\wedge\eta\big)\big(\tfrac{\partial}{\partial z_1},\tfrac{\partial}{\partial \overline{z_1}}\big)=\big(\tfrac{1}{4R}+m\big)i|z_2|^2$. 
Now since $dy_2\wedge dy_3=\tfrac{i}{2}\ddbar\big(|z_1z_2|^2\big)$, $dy_2\wedge dy_3\big(\tfrac{\partial}{\partial z_1},\tfrac{\partial}{\partial \overline{z_1}}\big)=\tfrac{i}{2}|z_2|^2$, thus $\zeta=\tfrac{1+4mR}{2R}=V$, and the lemma is proved. \cqfd 

~

One easily checks that $\eta$ is a connection 1-form away from 0 for the fibration given by $\varpi$; 
it is $S^1$-invariant (since $R$, $u$, $v$, $\tfrac{dz_1}{z_1}$, etc. are so), and at any point $p$ but $0\in\C^2$, since $\{\eta, dy_1, dy_2,dy_3\}$ is a basis of $T_p^*\C^2$ and $\varpi=(y_1,y_2,y_3)$, we necessarily have $T_p\C^2=\ker \eta+\ker T\pi$. 

Finally, the differential of $\eta$ has the expected shape:
 \begin{lem} \label{lem_deta}
  The differential of $\eta$ is given on $\C^2\backslash\{0\}$ by:
   \begin{equation*}
    d\eta=*_{\R^3}dV.
   \end{equation*}
 \end{lem}
\prf. Our 1-form $\eta$ is $S^1$-invariant and $\eta(\xi)$ is constant; by Cartan's formula, $0=\mathcal{L}_{\xi}\eta=\iota_{\xi}d\eta+d(\iota_{\xi}\eta)=\iota_{\xi}d\eta$, i.e.: the components of $d\eta$ in direction of the $dy_j\wedge\eta$ vanish. 
Moreover $d\omega_{\f}=0$ thus according to the latter proposition, $d\eta=\tfrac{\partial V}{\partial y_1}dy_2\wedge dy_3+\alpha_2 dy_3\wedge dy_1+\alpha_3 dy_1\wedge dy_2$. 
For the computation of $\alpha_2$ and $\alpha_3$, we do not have such a direct method. 
However, in order to simplify an explicit computation of $d\eta$, we can observe the following:
 \begin{equation*}
  4\eta=\Big(1+\frac{y_1}{R}\Big)d^c\log\big(|z_1|^2\big)-\Big(1-\frac{y_1}{R}\Big)d^c\log\big(|z_2|^2\big),
 \end{equation*}
by noticing simply that $u^2=R+y_1$ and $v^2=R-y_1$. 

Since $\log(|z_1|^2)$ and $\log(|z_2|^2)$ are pluriharmonic outside of $\{z_1=0\}$ and $\{z_2=0\}$, we thus have $d\eta=\tfrac{1}{4}d\big(\tfrac{y_1}{R}\big)\wedge d^c\log(|z_1z_2|^2)=\tfrac{1}{4}d\big(\tfrac{y_1}{R}\big)\wedge d^c\log(y_2^2+y_3^2)$. 
Now
 \begin{equation*}
  d\Big(\frac{y_1}{R}\Big)=\frac{1}{R^3}\Big((y_2^2+y_3^2)dy_1-y_1y_2dy_2-y_1y_3dy_3\Big)
 \end{equation*}
and 
 \begin{equation*}
  d^c\log(y_2^2+y_3^2)=I_1 d\log(y_2^2+y_3^2)=2\frac{y_2dy_3-y_3dy_2}{y_2^2+y_3^2}
 \end{equation*}
hence, after computations, $d\eta=-\tfrac{y_1}{2R^3}dy_2\wedge dy_3-\tfrac{y_2}{2R^3}dy_3\wedge dy_1-\tfrac{y_3}{2R^3}dy_1\wedge dy_2=\tfrac{\partial V}{\partial y_1}dy_2\wedge dy_3+\tfrac{\partial V}{\partial y_2}dy_3\wedge dy_1+\tfrac{\partial V}{\partial y_3}dy_1\wedge dy_2$. 

The lemma is proved, at least outside of $\{z_1z_2=0\}$, and there is no problem in extending the formula to $\C^2\backslash\{0\}$ by continuity. \cqfd

\subsection{Comparison of the Euclidean and the Taub-NUT metrics}  \label{sect_comp_ef}

 \subsubsection{Mutual control}

Euclidean and Taub-NUT geometries on $\C^2$ are radically different, in the sense that $\e$ and $\f$ are far from being globally mutually bounded. 
Indeed, euclidean balls volume grow proportionally to $r^4$, whereas after noticing that the function $R$ plays the role of the distance to 0 on $(\C^2,\f)$, 
we see that the volume of a big ball of radius $R$ of $(\C^2,\f)$ is essentially proportional to $R^3$. 

Another example of the geometric gap is given by the length of the circles of the $S^1$-action we used to check that $\f$ was the Taub-NUT metric; 
a circle passing through some $x\in\C^2\backslash\{0\}$ has length $2\pi r(x)$ when measured by $\e$, and length $2\pi V(x)^{-1/2}$ when measured by $\f$; 
this length tends to $\pi\sqrt{2/m}$ when $R(x)$ goes $\infty$ -- which gives us a geometric interpretation of the parameter $m$. 
We can nonetheless still compare the two metrics as follows:
 \begin{prop} \label{prop_comp_ef} 
  There exists some constant $C>0$ such that on $\C^2$ minus its unit ball, 
   \begin{equation*}
    C^{-1}r^{-2}\e \leq \f\leq Cr^2 \e.
   \end{equation*}
  \end{prop}
\prf. Since $\f=V(dy_1^2+dy_2^2+dy_3^2)+\tfrac{1}{V}\eta^2$, with $\eta=I_1Vdy_1$ and $dy_3=I_1dy_2$, we shall first evaluate $|dy_1|_{\e}$ and $|dy_2|_{\e}$. 
The latter one is immediate, from the identity $dy_2=\tfrac{i}{2}(z_1dz_2+z_2dz_1-\overline{z_1}d\overline{z_2}-\overline{z_2}d\overline{z_1})$, which provides $|dy_2|_{\e}=cr$. 
Now, we rearrange formulas \eqref{eqn_formula_uv} to write
 \begin{equation*}  \label{eqn_dy1}
  dy_1  = \frac{1}{2(1+4mR)}\big(e^{-4my_1}(\overline{z_1}dz_1+z_1d\overline{z_1})-e^{4my_1}(\overline{z_2}dz_2+z_2d\overline{z_2})\big). 
 \end{equation*}
This provides $|dy_1|_{\e}^2=\frac{c}{(1+4mR)^2}(|z_1|^2e^{-8my_1}+|z_2|^2e^{8my_1})$. 
But $|z_1|^2e^{-4my_1}=u^2$ and $|z_2|^2e^{4my_1}=v^2$, so $|dy_1|_{\e}^2=\frac{c}{(1+4mR)^2}(e^{-4my_1}u^2+e^{4my_1}v^2)=\frac{c'}{(1+4mR)^2}(R\cosh(4my_1)-y_1\sinh(4my_1))$. 
Now $R\cosh(4my_1)-y_1\sinh(4my_1)\leq R\cosh(4my_1)+y_1\sinh(4my_1)$ is obvious, and rearranging formulas \eqref{eqn_formula_uv} gives also
 \begin{equation}   \label{eq_formula_rRy1}
  2\big(R\cosh(4my_1)+y_1\sinh(4my_1)\big)=r^2,
 \end{equation}
so finally $|dy_1|^2_{\e}\leq c\tfrac{r^2}{R^2}$. 
Those estimates give us the bound $\f\leq Cr^2 \e$. 

On the other hand, at any point $p$, one can choose two 1-forms $e_1$ and $e_2$ so that $(e_1,Ie_1,e_2,Ie_2)$ is a basis of $T_p\C^2$, and that
 \begin{equation*}
  \e=e_1^2+(Ie_1)^2+e_2^2+(Ie_2)^2 \quad\text{and} \quad \f=\lambda\big(e_1^2+(Ie_1)^2\big)+\mu\big(e_2^2+(Ie_2)^2\big),
 \end{equation*}
$\lambda$, $\mu>0$. 
Since $\lambda\mu=\det^{\e}(\f)=1$, one has $\min(\lambda,\mu)=\max(\lambda,\mu)^{-1}$ in an obvious way. 
This allows us to convert the bound $\f\leq Cr^2 \e$ into a bound $\e\leq Cr^2 \f$. 
\cqfd

\begin{rmk}
 One can notice that with $R$ fixed, the function $y\mapsto R\cosh(4my)+y\sinh(4my)$ is convex and even on $[-R,R]$, so that its minimum is $R$, and its maximum is $Re^{mR}$. 
Consequently, we get from formula \eqref{eq_formula_rRy1} the inequality on $\C^2$:
 \begin{equation}  \label{eqn_rR}
  2R \leq r^2 \leq 2Re^{4mR},
 \end{equation}
which implies, among others, that $R$ is \textit{proper} on $\C^2$. 
This latter inequality is sharp, in the sense that $2R=r^2$ on $\{y_1=0\}=\{|z_1|=|z_2|\}$, and $2Re^{4mR}=r^2$ on $\{y_1=\pm R\}=\{z_1z_2=0\}$. 
\end{rmk}

  \subsubsection{Expressing euclidean objects in Taub-NUT vocabulary}

We give here some extra entries of the dictionary between the euclidean and the Taub-NUT settings on $\C^2$. 
More precisely, we give in Lemma \ref{lem_zeta} a vector field $\zeta$ helping to complete the dual frame of $(V^{-1/2}\eta,V^{1/2}dy_1,V^{1/2}dy_2,V^{1/2}dy_3)$ for $\f$. 
Then in Lemma \ref{lem_f_to_e}, we express the canonical frames of 1-forms and vector fields of $\e$, i.e. the $dx_j$ and the $\tfrac{\partial}{\partial x_j}$, in terms of those of $\f$. 

The essential point in those expressions lies in their computational consequences; indeed, they allow to compute objects like $\nabla^{\f}dx_j$, and estimate quantities like $\big|\nabla^{\f}dx_j\big|_{\f}$. 
Such estimates will be needed in next part, when constructing rough models of Ricci-flat Kähler metrics on the ALE spaces $X$ alluded to in Introduction. 

In paragraph \ref{prgrph_taubnut_C2}, we used the vector field $\xi$ on $\C^2$, which verified $\eta(\xi)=1$, $dy_j(\xi)=0$, $j=1,2,3$, and $dy_1(I_1\xi)=-\tfrac{1}{V}$, $\eta(I_1\xi)=dy_2(I_1\xi)=dy_3(I_1\xi)=0$. 
We shall complete our dual frame with the help of another vector field:
 \begin{lem} \label{lem_zeta}
  Define on $\C^2\backslash\{0\}$ the vector field
   \begin{equation*} \label{eq_df_zeta}
    \zeta=\frac{1}{2iR}\Big(e^{4my_1}\big(z_2\frac{\partial }{\partial \overline{z_1}}-\overline{z_2}\frac{\partial }{\partial z_1}\big)+e^{-4my_1}\big(z_1\frac{\partial }{\partial \overline{z_2}}-\overline{z_1}\frac{\partial }{\partial z_2}\big)\Big). 
   \end{equation*}
  Then $dy_2(\zeta)=1$ whereas $\eta(\zeta)=dy_1(\zeta)=dy_3(\zeta)=0$, and $dy_3(I_1\zeta)=1$ whereas $\eta(I_1\zeta)=dy_1(I_1\zeta)=dy_2(I_1\zeta)=0$. 
  Moreover, $[\xi,\zeta]=0$. 
 \end{lem}
\prf. We only need to check the first list of equalities, since $dy_3=I_1dy_2$ and $\eta=I_1Vdy_1$. 
Let us start with $dy_2(\zeta)=1$. 
Since $dy_2=\tfrac{1}{2i}(z_1dz_2+z_2dz_1-\overline{z_1}d\overline{z_2}-\overline{z_2}d\overline{z_1})$, we get 
 \begin{equation*}
  dy_2(\zeta)= \frac{1}{2R}\big(e^{4my_1}|z_2|^2+e^{-4my_1}|z_1|^2\big).
 \end{equation*} 
But $e^{4my_1}|z_2|^2=v^2$, $e^{-4my_1}|z_1|^2=u^2$, and $R=\tfrac{1}{2}(u^2+v^2)$, hence $dy_2(\zeta)=1$. 

To check that $dy_3(\zeta)=0$, use that $dy_3=-\tfrac{1}{2}(z_1dz_2+z_2dz_1+\overline{z_1}d\overline{z_2}+\overline{z_2}d\overline{z_1})$; this readily gives $dy_3(\zeta)=0$. 
Now for the equality $dy_1(\zeta)=0$, we use formula \eqref{eqn_dy1} to write $dy_1(\zeta)=\frac{1}{4iR(1+4mR)}(z_1z_2-\overline{z_1z_2}-z_1z_2+\overline{z_1z_2})=0$. 
The equality $\eta(\zeta)=0$ follows in the same way, using the formula for $\eta$ of Lemma \ref{lem_eta}. 

To see that $[\xi,\zeta]=0$, write $[\xi,\zeta]=\tfrac{d}{dt}e^{t\xi}\cdot\zeta\big|_{t=0}$ and recall that $e^{t\xi}\cdot(z_1,z_2)=(e^{it}z_1,e^{-it}z_2)$. 
Since $y_1$, $R$, as well as $\overline{z_2}\tfrac{\partial }{\partial z_1}$, $z_2\tfrac{\partial }{\partial \overline{z_1}}$, $\overline{z_1}\tfrac{\partial }{\partial z_2}$ and $z_1\tfrac{\partial }{\partial \overline{z_2}}$ are clearly invariant under this action, we get the result. 
\cqfd

 \begin{lem}  \label{lem_f_to_e}
  One has the following formulas for 1-forms:
   \begin{equation*}
    \begin{aligned}
     dx_1 &= Vx_1dy_1-x_2\eta  + \frac{e^{4my_1}}{2R}(x_4dy_2-x_3dy_3), \\
     dx_2 &= Vx_2dy_1+x_1\eta  + \frac{e^{4my_1}}{2R}(x_3dy_2+x_4dy_3), \\
     dx_3 &= -Vx_3dy_1+x_4\eta + \frac{e^{-4my_1}}{2R}(x_2dy_2-x_1dy_3), \\
     dx_4 &= -Vx_4dy_1-x_3\eta + \frac{e^{-4my_1}}{2R}(x_1dy_2+x_2dy_3) ;
    \end{aligned}
   \end{equation*}
  and for vector fields:
   \begin{equation*}
    \begin{aligned}
     \frac{\partial}{\partial x_1} &= -\frac{e^{-4my_1}}{2R}(x_2\xi+x_1I\xi)+(x_4\zeta-x_3I\zeta),  \\
     \frac{\partial}{\partial x_2} &= \frac{e^{-4my_1}}{2R}(x_1\xi-x_2I\xi) +(x_3\zeta+x_4I\zeta), \\
     \frac{\partial}{\partial x_3} &= \frac{e^{4my_1}}{2R}(x_4\xi+x_3I\xi)  +(x_2\zeta-x_1I\zeta), \\
     \frac{\partial}{\partial x_4} &= \frac{e^{4my_1}}{2R}(-x_3\xi+x_4I\xi)  +(x_1\zeta+x_2I\zeta).
    \end{aligned}
   \end{equation*}
 \end{lem}
\prf. We shall only see how those formulas arise for $dx_1$ and $\frac{\partial}{\partial x_1}$. 
Indeed, one can find the other identities using that $dx_2=I_1dx_1$, $dx_3=\tau^*dx_1$, $dx_4=I_1dx_3$, etc., on the euclidean side, and $\tau^{*}y_j=-y_j$, $\tau^*\eta=-\eta$, $\tau^*\xi=-\xi$, $\tau^*\zeta=-d\zeta$, etc., on the Taub-NUT side.

Write $dx_1=\alpha dy_1+\beta\eta+\gamma dy_2+\delta dy_3$. 
According to the duality between $(\xi,-VI_1\xi, \zeta,I_1\zeta)$ and $(\eta,dy_1,dy_2,dy_3)$, $dx_1(\xi)=\beta$, $dx_1(I\xi)=-\tfrac{\alpha}{V}$, $dx_1(\zeta)=\gamma$ and $dx_1(I_1\zeta)=\delta$. 
On the other hand, $dx_1(\xi)=\tfrac{1}{2}i(z_1-\overline{z_1})=-x_2$, 
$dx_1(I\xi)=-\tfrac{1}{2}(z_1+\overline{z_1})=-x_1$, 
$dx_1(\zeta)= \tfrac{1}{2iR}(e^{4my_1}\tfrac{1}{2}(z_2-\overline{z_2}+)+0)=\tfrac{e^{4my_1}}{2R}x_4$ 
and $dx_1(I\zeta)=\tfrac{i}{2R}(e^{4my_1}\tfrac{i}{2}(z_2+\overline{z_2})+0)=-\tfrac{e^{4my_1}}{2R}x_3$, hence the result. 

Similarly, write $\frac{\partial}{\partial x_1}=\alpha\xi+\beta I\xi+\gamma\zeta+\delta I\zeta$. 
Then $\alpha=\eta(\frac{\partial}{\partial x_1})=-e^{-4my_1}\frac{x_2}{2R}$, 
$\beta=-Vdy_1(\frac{\partial}{\partial x_1})=-e^{-4my_1}\frac{x_1}{2R}$, 
$\gamma=dy_2(\frac{\partial}{\partial x_1})= x_4$ and 
$\delta=dy_3(\frac{\partial}{\partial x_1})=-x_3$, as wanted. 
\cqfd

 \begin{rmk}
  With the previous formulas, it is easy to observe that
   \begin{equation*}
    dx_1\wedge dx_3+dx_4\wedge dx_2 = dy_2\wedge \eta+V dy_3\wedge dy_1 
   \end{equation*}
   and
   \begin{equation*}
    dx_1\wedge dx_4+dx_2\wedge dx_3 = dy_3\wedge \eta+V dy_1\wedge dy_2.
   \end{equation*}
  The reason for these identities is the following: one can check that the almost complex structures $J_2$ and $J_3$ defined by
   \begin{equation*}
    \left\{\begin{aligned}
            J_2Vdy_2=\eta,   \\
            J_2dy_3 = dy_1,
           \end{aligned} \right.
     \qquad \text{and} \qquad 
    \left\{\begin{aligned}
            J_3Vdy_3=\eta,   \\
            J_3dy_1 = dy_2,
           \end{aligned} \right.
   \end{equation*}
  are genuine complex structures (this follows from the harmonicity of $\varpi_*V$), and setting $J_1=I_1$, that we recover with $(J_j)$ the usual hyperkähler structure for the Taub-NUT metric. 
 Since now $\Theta$ has constant norm, hence is parallel, for $\f$, we get moreover that $\f(J_2\cdot,\cdot)=\mathfrak{Re}(\Theta)$ and $\f(J_3\cdot,\cdot)=\mathfrak{Im}(\Theta)$. 
 But we already know that $\Theta=\e(I_2\cdot,\cdot)+i\e(I_3\cdot,\cdot)$. 

 Notice moreover that when $m=0$, $J_2=I_2$ and $J_3=I_3$. 
 \end{rmk}

 \subsubsection{Derivatives}  \label{prgrpgh_drvativ}
  Consider the orthonormal frame $(e_j)_{j=0,\dots,3}$ of vector fields given by
   \begin{equation*}  \label{eqn_ej_frame}
    (e_0,e_1,e_2, e_3)=\big(V^{1/2}\xi, -V^{1/2}I_1\xi, V^{-1/2}\zeta, V^{-1/2}I\zeta\big)
   \end{equation*}
  away from 0. 
  In what follows, we will have to estimate the $\nabla^{\f}_{e_i}e_j$. 
  The following lemma will make the computations easier:
   \begin{lem}  \label{lem_lie_brckts}
    One has $[e_0,e_i]=\tfrac{y_i}{4R^3V^{3/2}}e_0$ for $i=1,2,3$, and
     \begin{equation*}
      [e_i,e_j]=\frac{1}{4R^3V^{3/2}}(y_ie_j - y_je_i +2y_ke_0)
     \end{equation*}
    for any triple $(i,j,k)\in \mathcal{I}=\big\{(1,2,3),(2,3,1),(3,1,2)\big\}$. 
As a consequence, 
     \begin{equation*}
      \nabla^{\f}e_0 = \frac{1}{4R^3V^{3/2}}\sum_{(i,j,k)\in\mathcal{I}} e_i\otimes\big(y_ke_j^*-y_je_k^*-y_ie_0^*\big)
     \end{equation*}
with $(e_0^*,e_1^*,e_2^*, e_3^*)=\big(V^{-1/2}\eta, V^{1/2}dy_1, V^{1/2}dy_2, V^{1/2}dy_3\big)$.
   \end{lem}

\begin{rmk}  \label{rmk_nabla_ei}
 Since $e_0=J_1e_1=J_2e_2=J_3e_3$ and the $(J_j)$ are parallel, an immediate induction shows that $(\nabla^{\f})^{\ell}e_i=O(R^{-(1+\ell)})$ for every $i\in\{0,1,2,3\}$ and $\ell\geq1$. 
One easily recovers from these estimates that $(\nabla^{\f})^{\ell}\riem^{\f}=O(R^{-(3+\ell)})$, which is one of the explanation for the terminology "Asymptotically Locally Flat" for $\f$.
\end{rmk}

~

\prf \textit{of Lemma \ref{lem_lie_brckts}}. 
Once the statement on the Lie brackets is proved, the formula for $\nabla^{\f}e_0$ follows from Koszul formula for the Levi-Civita connection $\nabla^{\f}$ and the orthonormality of the frame $(e_i)$. 
Moreover, because of the symmetric roles of $e_1$, $e_2$, $e_3$, we shall only see how to compute $[e_0,e_1]$ and $[e_1,e_2]$. 
 \begin{itemize}
  \item $[e_0,e_1]$: this bracket is rather easy to compute. 
   Recall that $e_0=V^{1/2}\xi$, $e_1=-V^{1/2}I_1\xi$, and $\xi$ is holomorphic for $I_1$, so that $[\xi,I_1\xi]=0$. 
   Moreover since $V$ is invariant under the $S^1$-action $\alpha\cdot(z_1,z_2)=(e^{i\alpha}z_1, e^{-i\alpha}z_2)$, $\xi\cdot V=0$ and $(I\xi)\cdot V=-V^{-1}\tfrac{\partial V}{\partial y_1}$. 
   On the other hand, 
    \begin{align*}
     [e_0,e_1]= & -[V^{1/2}\xi,V^{1/2}I_1\xi]= -V^{1/2}(\xi\cdot V^{1/2})I_1\xi + V^{1/2}(I_1\xi)\cdot V^{1/2}\xi - V [\xi,I_1\xi]                                   \\
              = &   -  \frac{1}{2} V^{-1}\frac{\partial V}{\partial y_1} \xi \quad  \qquad \text{since} \quad \xi\cdot V=0 \quad\text{and} \quad [\xi,I_1\xi]=0 \\
              = & \frac{1}{2V^{3/2}} \frac{\partial V}{\partial R}\frac{\partial R}{\partial y_1} e_0,
    \end{align*}
   hence the result, as $\tfrac{\partial V}{\partial R}=-\tfrac{1}{2R^3}$ and $\tfrac{\partial R}{\partial y_1}=\tfrac{y_1}{R}$. 

  \item $[e_1,e_2]$: this bracket is more delicate than the previous one. 
   We start with 
    \begin{equation*}
     [e_1,e_2]= -[V^{1/2}I_1\xi,V^{-1/2}\zeta]= -V^{1/2}(I_1\xi)\cdot V^{-1/2}\zeta + V^{-1/2}(\zeta\cdot V^{1/2})I_1\xi - [I_1\xi,\zeta].
    \end{equation*}
   We already know $(I_1\xi)\cdot V$, and similarly, $\zeta\cdot V=\tfrac{\partial V}{\partial y_2}$. 
   We are thus left with the computation of $[I_1\xi,\zeta]=\mathcal{L}_{I_1\xi}\zeta=\tfrac{d}{dt}\big|_{t=0}(e^{tI_1\xi})^*\zeta$. 
   But $e^{tI_1\xi}(z_1,z_2)=(e^{-t}z_1, e^tz_2)$; moreover, $\mathcal{L}_{I_1\xi}y_1=dy_1(I\xi)=-V^{-1}$, and consequently $\mathcal{L}_{I_1\xi}R=-\tfrac{y_1}{RV}$. 
   From the explicit formula \eqref{eq_df_zeta}, a direct computation yields
    \begin{equation*}
     \mathcal{L}_{I_1\xi}\zeta = \frac{|z_1|^2}{iR^2(1+4mR)}\big(z_2\frac{\partial }{\partial \overline{z_1}}-\overline{z_2}\frac{\partial }{\partial z_1}\big)
                               -\frac{|z_2|^2}{iR^2(1+4mR)}\big(z_1\frac{\partial }{\partial \overline{z_2}}-\overline{z_1}\frac{\partial }{\partial z_2}\big);
    \end{equation*}
   using formulas \eqref{eqn_formula_uv}, formulas from Lemma \ref{lem_f_to_e} and the identities $y_2=\tfrac{1}{2i}(z_1z_2-\overline{z_1z_2})$, $y_3=-\tfrac{1}{2}(z_1z_2+\overline{z_1z_2})$, this provides:
    \begin{equation*}
     \mathcal{L}_{I_1\xi}\zeta = -\frac{1}{2R^3V}(y_3\xi+y_2I_1\xi).
    \end{equation*}
   We can finally conclude that:
    \begin{align*} \hspace{0.8cm}
     [e_1,e_2]= &-\frac{1}{2V}V^{-1}\frac{y_1}{2R^3}\zeta + \frac{1}{2V}(-\frac{y_2}{2R^3})I_1\xi + \frac{1}{2R^3V}(y_3\xi+y_2I_1\xi) & \\ 
              = &\frac{1}{4R^3V^{3/2}}(y_2e_1-y_1e_2+2y_3e_0).                                                              &\hspace{0.8cm}\square 
    \end{align*}
\end{itemize}

~

\section{Application to $\mathcal{D}_k$ ALE spaces (minimal resolutions)}  \label{part_appl}
As  the Taub-NUT metrics of the previous part are $\mathcal{D}_k$-invariant, and since we have a biholomorphism 
between $\C^2/\mathcal{D}_k$ and its minimal resolution $X$ near infinity, a natural idea is to push-forward the 
potential of the Taub-NUT metric to $X$, compute the resulting (1,1)-form, and correct it into a global Ricci-flat Kähler form on $X$. 
We detail this program in the first section of this part; the following four sections are devoted to analytic considerations involved by the program.

  \subsection{Description of the program} \label{sctn_prgm}

  In order to construct ALF metrics on ALE gravitational instantons of dihedral type, the study of the Taub-NUT metric on $\C^2$ suggests the following program; here the positive parameter $m$ is fixed:
  \begin{enumerate}
   \item push the Taub-NUT metric $\f$ forward on a neighbourhood of infinity to $X$; if one assumes moreover that $X$ is the minimal resolution of $\C^2/\mathcal{D}_k$, one can identify $I_1$ and $I_1^X$, and then the pushed-forward metric is still Kähler where defined, 
    and one can take for it the push-forward of $\varphi$ as a potential.
   \item glue this metric at infinity with the ALE metric to get a Kähler metric on the whole $X$ (still with respect to $I_1^X$), with moreover a small Ricci form at infinity.
   \item correct the Ricci form on $X$ to get a Ricci-flat metric with ALF type asymptotics, and conclude by saying it is hyperkähler.
  \end{enumerate}

In a nutshell, if $(X,g)$ denotes an ALE gravitational instanton as above (minimal resolution) with volume form $\Omega_g$, and giving back its freedom to $m\in(0,+\infty)$, the result we get following this program is: 
 \begin{thm}  \label{thm_gRF}
  There exists a one-parameter family of hyperkähler metrics $g_{RF,m}$ on $X$, with volume form $\Omega_g$, 
  which are Kähler for $I_1^X$, indexed by a "mass" parameter $m\in(0,+\infty)$. 
  When $m$ is fixed, $(g_{RF,m}-\tilde{\f}_m)$ and $\nabla^{\tilde{\f}_m}(g_{RF,m}-\tilde{\f}_m)$ are $O(\tilde{R}^{-2-\delta})$ 
  for any $\delta\in (0,1)$. 
  If one denotes by $\omega_{RF,m}$ the Kähler form $g_{RF,m}(I_1^X\cdot,\cdot)$, 
  then $\omega_{RF,m}$ and $\omega_g$ are in the same cohomology class. 
 \end{thm}
In this statement, $\tilde{\f}_m$ is a smooth extension of the pull-back of $\pi_{*}\f$, with $\f$ computed with $m$, 
via the $I_1^X$-holomorphic identification between infinities of $\C^2/\mathcal{D}_k$ and $X$; $\tilde{R}$ is defined likewise.

Let us come back to the program itself; as suggested, we run it once $m$ is fixed, 
so we omit this parameter if not needed. 
Point 1. is rather clear, and does not need further detailed explanations. 
The gluing (point 2.) will be performed in next section. 
The correction of the Ricci-form (point 3.), which requires several steps of extra analytic work, is finally dealt with in sections \ref{sctn_RFnrinfty}, \ref{sctn_gge} and \ref{sctn_ccl_pgrm}.  

Nonetheless, we shall explain first  why we restrict ourselves to minimal resolutions from point 1, and do not deal with the general case of \textit{deformations}. 
Assuming we are in this general situation, call $\pi$ the projection $\C^2\to\C^2/\mathcal{D}_k$ and $\Phi$ an asymptotically tri-holomorphic diffeomorphism identifying $\C^2/\mathcal{D}_k$ and $X$ outside of compact subsets. 
We thus have three complex structures $I_j^X$, verifying the quaternionic relations, such that 
 \begin{equation*}
  \big|\nabla^{\ell}_{\e}\big(I_j-(\pi\circ\Phi)^*I^X_j\big)\big|_{\e}=O(r^{-4-{\ell}}),
 \end{equation*}
and this is sharp in general. 
In order to make a rough prototype of a Ricci-flat Kähler metric near infinity on $X$, one would like to take 
 \begin{equation*}
  (dI_1^Xd\varphi)(\cdot,I_1^X\cdot),
 \end{equation*}
where the $\varphi$ stands for $(\pi\circ\Phi)_*\varphi$. 

Call $\iota_1^X$ the difference $I_1-(\pi\circ\Phi)^*I^X_1$; we can think that for instance the estimate $|\iota_1^X|_{\e}=O(r^{-4})$ is rather nice. 
On the other hand, the error term between $(\pi\circ\Phi)_*\f$ and one can interpret $dI^X_1d(\pi\circ\Phi)^*\varphi$ on $\C^2$ as 
 \begin{equation*}
  (d\iota_1^Xd\varphi)(\cdot,I_1\cdot)+(dI_1d\varphi)(\cdot,\iota_1^X\cdot)+(d\iota_1^Xd\varphi)(\cdot,\iota_1^X\cdot).
 \end{equation*}

Now, the potential $\varphi$, which can be rewritten as $\varphi=\tfrac{1}{2}\big(R+m(R^2+y_1^2)\big)$, is thus $O(R^2)$, and easily $|d\varphi|_{\f}=O(R)$, $|\nabla^{\f}d\varphi|_{\f}=O(1)$. 
But because of the rather loose mutual control we have between $\e$ and $\f$ (Proposition \ref{prop_comp_ef}), this estimate converts into $|\iota_1^X|_{\f}=O(r^{-2})$ from inequality \eqref{eqn_rR}. 
At that point, we hence have $\big|\iota_1^Xd\varphi\big|_{\f}=O(r^{-2}R)$, which is  $O(1)$. 
What is to be estimated is however $\big|\nabla^{\f}(\iota_1^Xd\varphi)\big|_{\f}$, which is controlled by $|\iota_1^X|_{\f}\big|\nabla^{\f}d\varphi\big|_{\f}+|d\varphi|_{\f}\big|\nabla^{\f}\iota_1^X\big|_{\f}$. 

But the best we can do so far is : $|\iota_1^X|_{\f}\big|\nabla^{\f}d\varphi\big|_{\f}=O(r^{-2})=O(R^{-1})$, which may be good enough to go on, if there was not the bad estimate 
 \begin{equation} \label{eqn_nabla_iota}
  \big|\nabla^{\f}\iota_1^X\big|_{\f}=O(R^{-1})
 \end{equation}
giving $|d\varphi|_{\f}\big|\nabla^{\f}\iota_1^X\big|_{\f}=O(1)$.
This schematic estimation says that $(d\iota_1^Xd\varphi)(\cdot,I_1\cdot)$ is an error term with the same size as $\f$; 
there is no reason that the other error terms would compensate it, since 
 \begin{equation*}
  \big|(dI_1d\varphi)(\cdot,\iota_1^X\cdot)\big|_{\f}=\big|\f(I_1\cdot,\iota_1^X\cdot)\big|_{\f}\lesssim |\iota_1^X|_{\f} =O(R^{-1}),
 \end{equation*}
and
 \begin{equation*}
  \big|(d\iota_1^Xd\varphi)(\cdot,\iota_1^X\cdot)\big|_{\f}\lesssim |d\iota_1^Xd\varphi|_{\f}|\iota_1^X|_{\f}=O(R^{-1}).
 \end{equation*}

One could think though that we should get and $R^{-1}$ from the differentiation in \eqref{eqn_nabla_iota}; 
computations like those we make below (see e.g. the proof of Proposition \ref{prop_omegam}) however show that hope is not justified, and that the regularity we have on objects like $\iota_1^X$ for $\e$ does not propagate for such objects when considered via $\f$. 

These heuristic estimates are a bit rough; nonetheless, even a precise computation of the objects in play does not enable us to make them better. 
A way to bypass this difficulty may lie in the choice of the potential.

\subsection{Gluing the Kähler metrics}  \label{sect_gluing}

In this section, we work both on $X$ and on $\C^2$ for convenience. 
Since we fix the dihedral group $\mathcal{D}_k$ we work with, we fix the notation $\Gamma=\mathcal{D}_k$ for the rest of the article. 
To simplify the notations, we will not write the pull-backs explicitly when there is no risk of confusion. 

~

\paragraph{Volume forms on minimal resolutions.} 
On minimal resolutions $X$ of $\C^2/\Gamma$ for the pair $(I_1, I^X_1)$, endowed with their ALE metric $g$, one can assume that the volume form $\Omega_g$ is the same as the volume form $\Omega_{\e}$ induced by the euclidean metric on $X$ 
(this holds independently of the dihedral character of $\Gamma$); the volume form $\Omega_g$ in Theorem \ref{thm_gRF} is thus well-known.  
This also tells us that the $O(r^{-6})$ term in the volume form of the Eguchi-Hanson metric presented in Example \ref{expl_eh} is only an artifact due to the non-exhaustive expansion we give for the metric 
(nonetheless, this error term seems to be the correct error term in the more general case of deformations). 

~

\paragraph{The gluing.} In order to make a \textit{global} Kähler metric with $\f$ pushed-forward to $X$ minus a compact subset, we shall work with Kähler forms and more precisely with potentials, and use the already global $\omega_g$. 
We will also need a cut-off smooth non-decreasing function, $\chi$ say, such that
 \begin{equation} \label{eq_chi}
  \chi(t)=\left\{
           \begin{aligned}
            0 \quad &\text{if} \quad t\leq 0, \\
            1 \quad &\text{if} \quad t\geq 1.
           \end{aligned}
          \right.
 \end{equation}
We need a smooth nondecreasing convex function as well. 
We will denote by $\kappa:\R\to\R$ such a function, with the behaviour:
 \begin{equation*} \label{eq_kappa}
  \kappa(t)=\left\{
            \begin{aligned}
             0 \quad &\text{if} \quad t\leq 0, \\
             t \quad &\text{if} \quad t\geq 1
            \end{aligned}
           \right.
 \end{equation*}
(and as a consequence, one can choose $\chi=\kappa'$).

Moreover, we denote by $\varphi_0$ a smooth function on $X$ such that 
 \begin{equation*} \label{eq_df_alpha0}
  \alpha_0:=\omega_g-dd^c\varphi_0 =O(r^{-8}) \quad \text{and} \quad \big|(\nabla^g)^{\ell}(\omega_g-dd^c\varphi_0)\big|_g=O\big(r^{-8-\ell}\big) 
 \end{equation*}
for all $\ell\geq1$, where $r$ is a smooth positive function on $X$ extending the push-forward of the usual $r$ of $\C^2$. 
In particular, $dd^c\varphi_0\geq 0$ in the sense of $I_1^X$ (1,1)-forms near infinity.
Such a $\varphi_0$ can be chosen so that $\big|(\nabla^g)^{\ell}\varphi_0\big|_g=O\big(r^{2-\ell}\big)$ for all $\ell\geq0$. 

\begin{rmk}   \label{rmk_alpha0}
 Here one could have taken $\varphi_0$ so that the equality $\omega_{g}=dd^c\varphi_0$ exactly holds outside of a compact subset of $X$; this would have allowed us to avoid the analytic machinery we use in the following two sections. 
 We decided however to leave the error term $\alpha_0$, somehow artificially, because such a perturbation would have to be dealt with if performing our program in the case that $X$ is a deformation of $\C^2/\Gamma$; 
 we think indeed that the analytic considerations developed here would transpose nicely to that case, in order to get then a result generalizing Theorem \ref{thm_gRF} to any $\mathcal{D}_k$ ALE instanton. 
\end{rmk}

Recall that $\varphi$ is the Taub-NUT potential with fixed mass parameter $m$. 
We state:
 \begin{prop}  \label{prop_omegam}
  Take $K\geq0$ so that the identification $\Phi$ between infinities of $\C^2/\mathcal{D}_k$ and $X$ is defined on $\varphi\geq K$.
  Consider $r_0 \gg 1$, $\beta\in(0,1]$ and set 
   \begin{equation*}
    \Phi_m =  \kappa\circ(\varphi-K)- \chi\big((r-r_0)^{\beta}\big)\tilde{\varphi}_0,
   \end{equation*}
  where
   \begin{equation*}
    \tilde{\varphi}_0=\chi(r-r_0)\varphi_0
   \end{equation*} 
  Then if the parameters $K$ and $r_0$ (resp. $\beta$) are chosen big enough (resp. small enough), the symmetric 2-tensor $g_m$ associated to the $(1,1)$-form
   \begin{equation*}
    \omega_m:= \omega_g+dd^c\Phi_m
   \end{equation*}
  is well-defined on the whole $X$, is a Kähler metric, is ALF in the sense that 
   \begin{equation*}
    \big|(\nabla^{\f})^j(g_m-\f)\big|_{\f}=O(R^{-3}) \quad\text{for } j=0,1,2,3
   \end{equation*}
     and its volume form $\Omega_m$ verifies
   \begin{equation*}
    \big|(\nabla^{\f})^{\ell}(\Omega_m-\Omega_g)\big|_{\f}=O(R^{-3})
   \end{equation*}
for $\ell=0,1,2,3$, where $\Omega_g$ is the volume form of the ALE metric $g$.  
 \end{prop}

~

\noindent\prf. 
The proof splits into three steps, the first two corresponding to an analysis of each of the two components of the potential $\Phi_m$. 
In the first step, we deal with the component $\kappa\circ(\varphi-K)$ and adjust $K$; roughly speaking, we add an ALF part to an ALE Kähler metric, and get a metric which is asymptotically a sum ALE$+$ALF. 
This is the easy part. 
In the second step, we deal with $\chi\big((r-r_0)^{\beta}\big)\tilde{\varphi}_0$; this corresponds to killing the ALE part of the previous sum, in order to get an asymptotically ALF metric. 
Since we schematically subtract a metric to another one, this operation, where $r_0$ and $\beta$ are adjusted, needs to be carried out carefully to guarantee that the result is still a metric; this is why this second step is a bit more delicate than the first one. 
In the third step, we deal with the asymptotics of $g_m$ and of its volume form. 

~

\noindent\textit{First step.} 
Take $K$ as in the statement. 
A direct computation yields:
 \begin{equation*}
  dd^c\big(\kappa\circ(\varphi-K)\big)=\kappa''\circ(\varphi-K)d\varphi\wedge d^c\varphi+\kappa'\circ(\varphi-K)dd^c\varphi.
 \end{equation*}
Now since $\kappa$ is convex, $\kappa''\circ(\varphi-K)d\varphi\wedge d^c\varphi$ is always nonnegative in the sense of $I_1$ (or $I_1^X$) (1,1)-forms. 
Moreover, $dd^c\varphi$ is nothing but $\omega_{\f}$, with $\omega_{\f}=\f(I_1^X\cdot,\cdot)$ where this makes sense, so that finally, as soon as $K$ is big enough, $dd^c\big(\kappa\circ(\varphi-K)\big)\geq0$ on the whole $X$, and
 \begin{equation*}
  dd^c\big(\kappa\circ(\varphi-K)\big) = \omega_{\f}
 \end{equation*}
near infinity. 

Consequently, $\omega_g+dd^c\big(\kappa\circ(\varphi-K)\big)\geq \omega_g$ on the whole $X$, and this (1,1)-form equals $\omega_g+\omega_{\f}$ near infinity, namely on $\{\varphi\geq K+1\}$. 
We fix $K$ once for all. 

~

\noindent
\textit{Second step.} 
Let us choose now $r_0$ so that $\{r\geq r_0\}$ makes sense on $X$ and is a subset of $\{\varphi\geq K\}$. 
One has, if setting $\tilde{\alpha}_0:=\omega_g-dd^c\tilde{\varphi}_0$, which equals $\alpha_0$ defined in \eqref{eq_df_alpha0} on $\{r\geq r_0+1\}$:
 \begin{equation*}
  \begin{aligned}
   \omega_g-dd^c_X\big(\chi\big((r-r_0)^{\beta}\big)\tilde{\varphi}_0\big)  =& \chi\big((r-r_0)^{\beta}\big) \tilde{\alpha}_0 
                                                                            + \big(1-\chi\big((r-r_0)^{\beta}\big)\big)\omega_g      \\
                                                                          &- \beta\chi'\big((r-r_0)^{\beta}\big)(r-r_0)^{\beta-1}\Xi_{\beta}  \\
                                                                          &- \beta^2 \chi''\big((r-r_0)^{\beta}\big)(r-r_0)^{2\beta-2}\Upsilon
  \end{aligned}
 \end{equation*}
where
 \begin{equation*}
  \left\{\begin{aligned}
          \Xi_{\beta}=& d\tilde{\varphi}_0\wedge d^cr+ dr\wedge d^c_X\tilde{\varphi}_0 
                      +\tilde{\varphi}_0 dd^cr+(\beta-1)(r-r_0)^{-1}\tilde{\varphi}_0dr\wedge d^cr   \\
          \Upsilon   =& \tilde{\varphi}_0dr\wedge d^cr 
         \end{aligned}
  \right.
 \end{equation*}
We deal separately with the different summands: 
 \begin{itemize}
  \item[$\bullet$] $\big(1-\chi\big((r-r_0)^{\beta}\big)\big)\omega_g$ is everywhere nonnegative; 
  \item[$\bullet$] $\chi\big((r-r_0)^{\beta}\big) \tilde{\alpha}_0$ is $O(r^{-6})$ for $g$ (or $\e$) (independently of $\beta$) and is zero on $\{r\leq r_0\}$, so is bounded below by $-\tfrac{1}{2}\omega_{\f}$ on $\{r\geq r_0\}$ provided that $r_0$ is fixed big enough; 
   we fix such an $r_0$ once for all; 
  \item by construction, there exists some constant $C$ (depending on $r_0$, which is fixed, but independent of $\beta$), such that:
   \begin{equation*}
    \big|\tilde{\varphi}_0\big|\leq C(r-r_0)^2,\quad \big|d\tilde{\varphi}_0\big|_{g}\leq C(r-r_0), \quad \big|dd^c\tilde{\varphi}_0\big|_{g}\leq C;
   \end{equation*}
taking a bigger constant (still independent of $\beta$) if necessary to take into account the contributions of $dr=O(1)$, $dd^cr=O(\tfrac{1}{r})$, $\chi'\big((r-r_0)^{\beta}\big)$ and so on, one thus has:
 \begin{align*}
  \Big|\beta\chi'\big((r-r_0)^{\beta}\big)(r-r_0)^{\beta-1}\Xi_{\beta}
       + \beta^2 \chi''\big((r-&r_0)^{\beta}\big)(r-r_0)^{2\beta-2}\Upsilon \Big|_{g}   \\
                               &  \leq C(\beta(r-r_0)^{\beta}+\beta^2(r-r_0)^{2\beta})
 \end{align*}
on the annulus where the left-hand side member may not vanish identically, i.e. $\{r_0\leq r\leq r_0+1\}$. 
Now on this region $(r-r_0)^{\beta}\leq1$, so the latter bound becomes $C(\beta+\beta^2)$. 

Take $c>0$ small enough so that $\omega_{\f}\geq \tfrac{c}{r}\omega_g$ on $\{r\geq r_0\}$, and choose $\beta$ small enough so that $\beta+\beta^2\leq \frac{c}{4C(r_0+1)^2}$. 
Then 
 \begin{equation*}
  \beta\chi'\big((r-r_0)^{\beta}\big)(r-r_0)^{\beta-1}\Xi_{\beta} + \beta^2 \chi''\big((r-r_0)^{\beta}\big)(r-r_0)^{2\beta-2}\Upsilon    
                                                                                                                           \geq -\frac{1}{4}\omega_{\f}
 \end{equation*} 
on $\{r_0\leq r\leq r_0+1\}$, and the left-hand side member vanishes elsewhere. 
 \end{itemize} 
 
Finally, on $X\backslash \{r>r_0\}$, $\omega_m:= \omega_g+dd^c\Phi_m\geq \omega_g$; 
on the annulus $\{r_0\leq r\leq r_0+1\}$, $\omega_m\geq \tfrac{1}{4}\omega_{\f}$; on the neighbourhood of infinity $\{r\geq r_0+1\}$, $\omega_m\geq \tfrac{1}{2}\omega_{\f}$. 
In other words, $\omega_m$ is positive, hence Kähler, on the whole $X$.

~

\noindent
\textit{Third step.} 
Near infinity, namely on $\{r\geq r_0+1\}$, $\omega_m= dd^c\varphi+\alpha_0=\omega_{\f}+\alpha_0$. 
We are thus left with proving that $\big|(\nabla^{\f})^{\ell}\alpha_0\big|_{\f}=O(R^{-3})$ for $\ell\in\{0,\dots,3\}$. 

The case $\ell=0$ does not require much work: write $|\alpha_0|_{\f}\leq |g|_{\f}|\alpha_0|_{g}$ ($\alpha_0$ is a $(0,2)$ tensor) and $|g|_{\f}=O(r^2)$, since $g\sim\e$, so $|\alpha_0|_{\f}=O(r^{-6})$, which is a $O(R^{-3})$. 

In order to evaluate $\nabla^{\f}\alpha_0$, let us work on $\C^2\backslash B$, where we use global coordinates adapted to $\e$, and also of global coordinates and connection $1$-from adapted to $\f$. 
We still denote the pull-back of $\alpha_0$ by $\alpha_0$. 
Let us write on $\C^2\backslash B$
 \begin{equation*}
  \alpha_0 = \sum_{j,k}\alpha_{jk}dx_j\wedge dx_k
 \end{equation*}
so that 
 \begin{equation} \label{eqn_nabalf_alpha0}
  \nabla^{\f}\alpha_0= \sum_{j,k}d\alpha_{jk}\otimes(dx_j\wedge dx_k)+ \sum_{j,k}\alpha_{jk}\nabla^{\f}(dx_j\wedge dx_k).
 \end{equation}
On the one hand, $\nabla^{\e}\alpha_0=d\alpha_{jk}\otimes(dx_j\wedge dx_k)$ is a $O(r^{-9})$ for $\e$, i.e. $d\alpha_{jk}$ are so. 
Since they are 1-forms, one gets that they are $O(r^{-8})$ for $\f$ (we lose the $r^{-1}$ we won by differentiating $\alpha_{jk}$), and then that $d\alpha_{jk}\otimes(dx_j\wedge dx_k)$ is a $O(r^{-6})$, hence a $O(R^{-3})$, for $\f$. 
On the other hand, we have to estimate the $\nabla^{\f}dx_j$; we state: 
 \begin{lem} \label{lem_nablafdxj}
  For $j\in\{1,\dots,4\}$, one has $|\nabla^{\f}dx_j|_{\f}=O(r)$. 
  More generally, for any $\ell\geq2$, one has $\big|(\nabla^{\f})^{\ell}dx_j\big|_{\f}=O(r)$.   
 \end{lem}
We show these estimates after the end of the current proof. 

~

For now, we conclude from these estimates and formula \eqref{eqn_nabalf_alpha0} that $|\nabla^{\f}\alpha_0|_{\f}=O(r^{-6})$, which is a $O(R^{-3})$, and we transpose this estimation on $X$ near its infinity. 
It is an easy induction to prove that $\big|(\nabla^{\f})^{\ell}\alpha_0\big|_{\f}=O(R^{-3})$ for $\ell=2,3$. 

Let us deal finally with the final on the volume form of $\omega_m$. 
One writes
 \begin{equation*}
  \Omega_m-\Omega_g=(\Omega_m-\Omega_{\f})+(\Omega_{\f}-\Omega_{\e})+(\Omega_{\e}-\Omega_g)
 \end{equation*}
The asymptotics we proved on $\omega_{m}-\omega_{\f}$, hence on $g_m-\f$, provide $\big|(\nabla^{\f})^{\ell}(\Omega_m-\Omega_g)\big|_{\f}=O(R^{-3})$ for $\ell=0,\dots,3$. 
Moreover, the central form $\Omega_{\f}-\Omega_{\e}$ vanishes, and thanks to our assumption that $X$ is a minimal resolution of $\C^2/\Gamma$, so does $\Omega_{\e}-\Omega_g$. 
\cqfd

~

\noindent
\textit{Proof of Lemma \ref{lem_nablafdxj}}. 
We shall only deal with the estimate on $\nabla^{\f}dx_1$, since it will be rather clear from our computations that those of the $(\nabla^{\f})^{\ell}dx_1$, $\ell\geq2$, would be established in a similar way. 
The estimates on the $(\nabla^{\f})^{\ell}dx_j$, $j=2,3,4$, would furthermore be obtained by symmetry. 

From Lemma \ref{lem_f_to_e}, using the orthonormal frame $(e_j^*)$ given by in Lemma \ref{lem_lie_brckts}, one has
 \begin{equation*}
  dx_1 = -x_2V^{1/2}e_0^* +V^{1/2}x_1e_1^*+\frac{e^{4my_1}}{2R}(V^{-1/2}x_4 e_2^*- V^{-1/2} x_3e_3^*).
 \end{equation*}
Thus
 \begin{equation} \label{eqn_nablafdx1}
  \begin{aligned}
   \nabla^{\f}dx_1 = &-dx_2\otimes V^{1/2}e_0^* +V^{1/2}dx_1\otimes e_1^*+\frac{e^{4my_1}}{2R}(V^{-1/2}dx_4\otimes e_2^*- V^{-1/2}dx_3\otimes e_3^*) \\
                    -&x_2\nabla^{\f}(V^{1/2}e_0^*)+x_1 \nabla^{\f}(V^{1/2}e_1^*) +\frac{e^{4my_1}}{2R}\big(x_4\nabla^{\f}(V^{-1/2} e_2^*)- x_3\nabla^{\f}(V^{-1/2} e_3^*)\big) \\
                    +&d\big(\frac{e^{4my_1}}{2R}\big)\otimes (V^{-1/2}x_4 e_2^*- V^{-1/2} x_3e_3^*). 
  \end{aligned}
 \end{equation}
One handles this formula in the following way (estimations are made w.r.t. $\f$):
 \begin{itemize}
  \item $dx_1,dx_2=O(r)$ so $-dx_2\otimes V^{1/2}e_0^* +V^{1/2}dx_1\otimes e_1^*=O(r)$; 
  \item $\tfrac{e^{4my_1}}{2R}dx_4=  -\frac{e^{4my_1}}{2R}(Vx_4dy_1+x_3\eta) + x_1dy_2+x_2dy_3$, so 
   \begin{align*}
    \big|\frac{e^{4my_1}}{2R}dx_4\big|_{\f}=&\Big(V|z_2|^2\big(\frac{e^{4my_1}}{2R}\big)^2+V^{-1}|z_1|^2\Big)^{1/2} 
                                           =  e^{2my_1} \big(V\frac{v^2}{4R^2}+V^{-1}u^2\big)^{1/2}                 \\
                                           =& O(r),
   \end{align*}
  since $e^{2my_1}=O(\tfrac{r}{R^{1/2}})$ from formula \eqref{eqn_rR}; similarly, $\tfrac{e^{4my_1}}{2R}dx_3=O(r)$. 
  Hence $\tfrac{e^{4my_1}}{2R}(V^{-1/2}dx_4\otimes e_2^*- V^{-1/2}dx_3\otimes e_3^*)=O(r)$, and thus the first line of the RHS of \eqref{eqn_nablafdx1} is a $O(r)$ ;
  \item $\nabla^{\f}(V^{1/2}e_0^*)$, $\nabla^{\f}(V^{1/2}e_1^*)$, $\nabla^{\f}(V^{-1/2} e_2^*)$ and $\nabla^{\f}(V^{-1/2} e_3^*)=O(1/R^2)$ from Lemma \ref{lem_lie_brckts}; 
  $x_1$, $x_2=O(r)$, and $\tfrac{e^{4my_1}}{2R}x_3$, $\tfrac{e^{4my_1}}{2R}x_4=O(\tfrac{r}{R})$, as $(e^{4my_1})^2|z_2|^2=e^{4my_1}v^2=O(r^2)$. 
  The second line of the RHS of \eqref{eqn_nablafdx1} is a $O(\tfrac{r}{R^2})$ ;
  \item $d\big(\frac{e^{4my_1}}{2R}\big)=\tfrac{mRdy_1-dR}{2R^2}e^{4my_1}$, from which we deduce as above that $d\big(\tfrac{e^{4my_1}}{2R}\big)\otimes (V^{-1/2}x_4 e_2^*- V^{-1/2} x_3e_3^*)=O(\tfrac{r}{R})$; 
  the last line of the RHS of \eqref{eqn_nablafdx1} is a $O(\tfrac{r}{R})$. 
 \end{itemize}
Collecting these estimates, we can conclude that $|\nabla^{\f}dx_1|_{\f}=O(r)$. \cqfd

\subsection{Making the metric Ricci-flat near infinity} \label{sctn_RFnrinfty}

We shall see now that we can correct our metric $g_m$ to a metric which is Ricci-flat in a neighbourhood of infinity. 
Indeed, before making the metric Ricci-flat \textit{on the whole }$X$, which is the final aim of point 4. in the program sketched in paragraph \ref{sctn_prgm}, we need to control our metric up to higher orders.  

The price to pay to correct $g_m$ is a slight loss of accuracy in the asymptotics of the metric for low orders, which may seem paradoxical since we want to improve regularity. 
However, a way to obtain the regularity we need, which is \textit{local near infinity}, once we have a  Ricci-flat metric outside a compact subset, is putting it in an adequate gauge, which corresponds to looking at it with a better chart at infinity. 
This compensates the slight loss alluded to above, and provides actually asymptotics at any order. 

For now, we settle the problem of getting Ricci-flatness outside of a compact subset of $X$; we perform the change of gauge in next section. 

The first point uses a standard implicit functions theorem; in order to apply it properly, we have to introduce quickly weighted Hölder spaces. 
Notice that thanks to the asymptotics of Proposition \ref{prop_omegam}, and since $\f$ has a nonnegative injectivity radius, so does the metric $g_m$ produced in this proposition; denote this injectivity radius by $\inj_{g_m}$. 

 \begin{df}[Weighted Hölder spaces.] \label{df_wtd_hldr_sp} 
  Let $(k,\alpha)\in\N\times[0,1)$, and $\delta\in\R$. 
  For any $C^{0,\alpha}$ tensor $u$, define
   \begin{equation} \label{eq_hldrmdls}
    [u]^{\alpha}_{\delta}(x)= 
                            \sup_{\substack{y\in X,\\ d_{g_m}(x,y)<\inj_{g_m}}}
                             \Big|\min(\tilde{R}(x)^{\alpha+\delta},\tilde{R}(y)^{\alpha+\delta})
                             \frac{u(x)-u(y)}{d_{g_m}(x,y)^{\alpha}}\Big|_{g_m}
   \end{equation}
  where $\tilde{R}$ is a positive smooth extension of $R$ on $X$. 

  We define the weighted Hölder space $C^{k,\alpha}_{\delta}$ as the space of $C^{k,\alpha}$ functions $f$ such that 
   \begin{equation*}
    \|f\|_{C^{k,\alpha}_{\delta}}:= \big\|\tilde{R}^{\delta}f\big\|_{C^0}+\cdots+\big\|\tilde{R}^{k+\delta}(\nabla^{g_m})^kf\big\|_{C^0}
                                   +\sup_X \big[\tilde{R}^{k+\delta}(\nabla^{g_m})^kf\big]^{k,\alpha}_{\delta}<+\infty
   \end{equation*}
  where the norms $\|\cdot\|_{C^0}$ are computed with $g_m$. 
  We endow this space with this $\|\cdot\|_{C^{k,\alpha}_{\delta}}$ norm. 
 \end{df}
  
In the definition of $[\cdot]^{\alpha}_{\delta}$ (formula \eqref{eq_hldrmdls}), the difference $\big(u(x)-u(y)\big)$ is interpreted in a usual way, i.e. bringing $u(y)$ to $x$ with the parallel transport along the minimizing geodesic joining $y$ to $x$. 

 \begin{rmk}
  Notice that in view of the asymptotics of $\omega_m$, one would have obtained the same $C^{k,\alpha}_{\delta}$ if defined with (a smooth extension of) $\f$, as long as $k+\alpha+\delta\leq 3$ ($\delta\geq 0$). 
 \end{rmk}

We recall that $g$ is the original hyperkähler ALE metric on $X$.

 \begin{prop} \label{prop_RFatinfty}
  Take $(\alpha_1,\delta_1)\in(0,1)^2$ such that $\alpha_1+\delta_1<1$. 
  Then there exists $R_0=R_0(\alpha_1,\delta_1)$ and a smooth function $\psi\in C^{2,\alpha}_{\delta}\cap C^{3,\alpha}_{\delta-1}\cap C^{3,\alpha}_{\delta-2}$ such that 
  $\omega_{\psi}:=\omega_m+dd^c\psi$ is positive and satisfies
   \begin{equation*} \label{eq_omgpsiRF}
    \frac{1}{2}\omega_{\psi}^2 = \Omega_g,
   \end{equation*}
  and is thus Ricci-flat, on $\{R\geq R_0+1\}$. 
 \end{prop}
 \prf. Consider the cut-off function $\chi$ of the previous section (see equations \eqref{eq_chi}), take $R_0\gg 1$, set $\chi_{R_0}=\chi(R-R_0)$ on $\{R\geq R_0\}$ and extend it by 0 in elsewhere; this makes sense provided $R_0$ is large enough.
  
 If we solve the problem
  \begin{equation}   \label{omgpsiRF_pb}
   \big(\omega_m+dd^c\psi\big)^2= (1-\chi_{R_0})\omega_m^2 + 2\chi_{R_0} \Omega_g
  \end{equation}
with $\psi\in C^{2,\alpha_1}_{\delta}\cap C^{3,\alpha_1}_{\delta_1-1}\cap C^{4,\alpha_1}_{\delta_1-2}$ where $\alpha_1,\delta_1>0$ and $\alpha_1+\delta_1<1$, such that $\omega_m+dd^c\psi>0$, 
then we are done up to the positivity assertion, since the smoothness of $\psi$, which is purely local, will follow at once from local ellipticity arguments. 

Now, we want to solve our problem \eqref{omgpsiRF_pb} with the implicit function theorem. 
Indeed, the linearisation of the Monge-Ampère operator
 \begin{equation*}
  \begin{aligned}
   C^{k+2,\beta}_{\nu} & \longrightarrow      C^{k,\beta}_{\nu+2}                                        \\
        \psi           & \longmapsto     \frac{\big(\omega_m+dd^c\psi\big)^2-\omega_m^2}{\omega_m^2}
  \end{aligned}
 \end{equation*}
($k\geq0$, $\beta\in(0,1)$, $\nu>0$) is, up to the sign, the Laplacian operator $\Delta_{g_m}$ of $g_m$. 
Moreover, we know that $g_m$ is $C^{0,\alpha}_{2+\delta}\cap C^{1,\alpha}_{1+\delta}\cap C^{2,\alpha}_{\delta}$ close to $\f$, with $\alpha$, $\delta\in(0,1)$, $\alpha+\delta\leq1$ and $\alpha>\alpha_1$, $\delta>\delta_1$. 
Let us quote the following lemma:
 \begin{lem} \label{lem_lapl}
  Take $k\in\N$, $\nu\in\R$ such that $k+\nu+\alpha_1<3$. 
  The operator $\Delta_{g_m}:C^{k+2,\alpha_1}_{\nu}\to C^{k,\alpha_1}_{\nu}$ is an isomorphism if $\nu\in(0,1)$, and surjective with kernel reduced to the constants if $\nu\in(-2,0)$. 
 \end{lem}
We postpone the proof of this lemma below. 

We get from this lemma isomorphisms $\Delta_{g_m}:  C^{2,\alpha_1}_{\delta_1}\to C^{0,\alpha_1}_{2+\delta_1}$, $E^{2,\alpha_1}_{\delta_1-1}(x_0)\to C^{0,\alpha_1}_{1+\delta_1}$ and $E^{3,\alpha_1}_{\delta_1-2}(x_0)\to C^{1,\alpha_1}_{\delta_1}$, 
where $E^{\bullet}_{*}(x_0)$ is the space of $C^{\bullet}_{*}$ functions vanishing at some fixed point $x_0\in X$.

So the last ingredient we need to apply the implicit function theorem to our Monge-Ampère operator is making the right-hand side member in \eqref{omgpsiRF_pb} arbitrarily close to $\omega_m^2$ in $C^{0,\alpha}_{2+\delta}$, $C^{1,\alpha}_{1+\delta}$ and $C^{2,\alpha}_{\delta}$ norms. 
In other terms, we want the $C^{0,\alpha_1}_{2+\delta_1}$, $C^{1,\alpha_1}_{1+\delta_1}$ and $C^{2,\alpha_1}_{\delta_1}$ norms of 
  \begin{equation*}
  \omega_m^2- \big((1-\chi_{R_0})\omega_m^2 + 2\chi_{R_0} \Omega_g\big)=2\chi_{R_0}(\Omega_m-\Omega_g)
 \end{equation*}
to be arbitrarily small. 
Since $(\nabla^{\f})^{\ell}(\Omega_m-\Omega_g)=O(R^{-3})$ and $\alpha_1+\delta_1<1$, an easy interpolation provides this estimate for $R_0$ big enough. 

We can thus apply simultaneously the implicit function theorem for our Monge-Ampère operator, 
and get $\psi_0\in C^{2,\alpha_1}_{\delta_1}$ $\psi_1\in \psi_0(x_0)+E^{3,\alpha_1}_{\delta_1-1}(x_0)\subset C^{2,\alpha_1}_{\delta_1-1}$ and $\psi_2\in \psi_0(x_0)+E^{4,\alpha_1}_{\delta_1-2}(x_0)\subset C^{3,\alpha_1}_{\delta_1-2}$ satisfying problem \eqref{omgpsiRF_pb}. 
Before we see that $\psi_0=\psi_1=\psi_2$, we shall check that the $\omega_{\psi_j}$ are positive. 

In any cases, $\omega_{\psi_j}=\omega_m+o(1)$, this $o(1)$ being computed with $\omega_m$, so that $\omega_{\psi_j}$ is positive near infinity. 
Now, its volume form, which we can interpret as $\det^{\omega_m}(\omega_{\psi_j})\omega_m^2$, never vanishes, and hence neither do its eigenvalues with respect to $\omega_m$. 
Thus the $\omega_{\psi_j}$ are indeed positive as $I_1^X$ (1,1)-forms on the whole $X$. 

To see that $\psi_1=\psi_2$, let us set $\psi_{12}=\psi_1-\psi_2$, call $g_{\psi_1}$ (resp. $g_{\psi_2}$) the metric associated to $\omega_{\psi_1}$ (resp. $\omega_{\psi_2}$), and observe that 
 \begin{equation*}
  (\omega_{\psi_1}+\omega_{\psi_2})\wedge dd^c_X\psi_{12}=(\omega_{\psi_1}+\omega_{\psi_2})\wedge (\omega_{\psi_1}-\omega_{\psi_2}) 
                                                         =  \omega_{\psi_1}^2-\omega_{\psi_2}^2 
                                                         =0.
 \end{equation*}
On the other hand, $\omega_{\psi_1}+\omega_{\psi_2}$ is a Kähler form, namely that associated to $g_{12}:=g_{\psi_1}+g_{\psi_2}$, so actually $(\omega_{\psi_1}+\omega_{\psi_2})\wedge dd^c_X\psi_{12}=-\Delta_{g_{12}}\psi_{12}\vol^{g_{12}}$. 
In other words, $\psi_{12}$ is harmonic for $g_{12}$. 
But $\psi_{12}=O(R^{2-\delta})$, and $g_{12}$ is close enough to the Kähler metric $2g_{\psi_1}$, which is itself close enough to $2g_m$ say, so that we know that $\psi_{12}$ is constant. 
Since $\psi_{12}(x_0)=0$, $\psi_{12}\equiv 0$, i.e. $\psi_1=\psi_2$. 
The identity $\psi_0=\psi_1$ is deduced in the same way. 
\cqfd

~

Notice however that the last argument in the latter proof works thanks to the complex dimension 2. 

~

\noindent
\textit{Proof of Lemma \ref{lem_lapl}}. 
The statement for $\Delta_{g_m}:C^{k+2,\alpha'}_{\nu}\to C^{k,\alpha'}_{\nu}$ with $\nu\in(0,1)$ is done in the appendix of \cite{bm}. 
Using the same techniques, and the facts that $0$ is the first critical weight, authorizing constants, and $-2$ is the second critical weight, we get the assertion on $\Delta_{g_m}:C^{k+2,\alpha'}_{\nu}\to C^{k,\alpha'}_{\nu}$ with $\nu\in(-2,0)$. 
The fact that the second critical weight is $-2$, and not $-1$, is due to the action of $\mathcal{D}_k$, which makes $(\C^2/\mathcal{D}_k,\f)$, hence $(X,\tilde{\f})$, into a circle fibration over $\R^3/{\pm}$ at infinity.
\cqfd

\subsection{Improvement of the asymptotics} \label{sctn_gge}

When one has a Ricci-flat metric, it is a somehow usual trick to put it in an adequate gauge to improve the regularity we can state on it. 
Roughly speaking, this process corresponds to look at the metric in a better chart, in which the artifacts of its construction disappear, letting one use inductive ellipticity arguments to reach the desired regularity. 

We hence want to look at our metric $g_{\psi}$ pulled-back by some diffeomorphism, or, which is equivalent for our purpose, to look at its pull-back on $\C^2$ pulled-back once more by a diffeomorphism of $\C^2$; we work on this space instead of $X$ for minor technical reasons. 
Before explaining what will guide us for the choice of such a diffeomorphism, let us specify the class we shall take it in.

 \subsubsection{The diffeomorphisms}
  
We recall that the orthonormal frames $(e_i)$ and $(e_i^*)$ are given in paragraph \ref{prgrpgh_drvativ}. 

  \begin{df}
   Let $(k,\alpha)\in\N^*\times(0,1)$, and let $\nu>0$. 
   We denote by $\diff^{k,\alpha}_{\nu}$ the class of diffeomorphisms $\phi$ of $\C^2$ such that:
    \begin{itemize}
     \item $\phi$ has regularity $(k,\alpha)$;
     \item there exists a constant $C$ such that for any $x\in \C^2$, $d_{\f}\big(x,\phi(x)\big)\leq C\big(1+R(x)\big)^{-\nu}$;
     \item take $R_0\geq 1$ such that for any $x\in\{R\geq R_0\}$, then $d_{\f}\big(x,\phi(x)\big)\leq \inj_{\f}$. 
      Denote by $\gamma_x:[0,1]\to\C^2$ the minimizing geodesic for $\f$ joining $\phi(x)$ to $x$ and by $p_{\gamma_x}$ the parallel transport along $\gamma_x$. 
      Consider the maps $\phi_{ij}: \{R\geq R_0\} \to \R$ given by 
       \begin{equation*}
        \phi_{ij}(x)=e_{i}^*\big((T_x\Phi\circ p_{\gamma_x}(1)-\id)(e_{j})\big),
       \end{equation*}
      and extend them smoothly to $\{R\leq R_0\}$. We then ask $\|\phi_{ij}\|_{C^{k-1,\alpha}_{\nu+1}}$ to be finite. 
    \end{itemize}
   We endow $\diff^{k,\alpha}_{\nu}$ with the obvious topology.
  \end{df}
In this definition the involved weighted Hölder spaces are the analogues for $\f$ of those of section \ref{sctn_RFnrinfty} for $g_m$. 

This definition gives a concrete idea of what is a diffeomorphism of class $C^{k,\alpha}_{\nu}$, and makes the following lemma quite intuitive:
 \begin{lem}
  For any $Y\in C^{k,\alpha}_{\nu}(\C^2)$ with $\|Y\|_{C^{k,\alpha}_{\nu}}$ small enough, the map
   \begin{equation*}
    \phi_Y : x\longmapsto \exp^{\f}_x\big(Y(x)\big)
   \end{equation*}
  is in $\diff^{k,\alpha}_{\nu}$. 
  Conversely, any $\phi\in \diff^{k,\alpha}_{\nu}$ close enough to the identity has the form $\phi_Y$
  for a unique $Y\in C^{k,\alpha}_{\nu}$. 

  In other words, one can parametrize a neighbourhood of identity in $\diff^{k,\alpha}_{\nu}$ by a neighbourhood of $0$ in $C^{k,\alpha}_{\nu}$. 
 \end{lem}
 \prf. This statement is standard, so we will content ourselves with a word on the injectivity of the $\phi_Y$. 
 Suppose $\|Y\|_{C^{k,\alpha}_{\nu}}\leq 1$. 
 Take a triplet $(y,x_1,x_2)$ of points of $X$ such that $y=\phi_Y(x_1)=\phi_Y(x_2)$. 
 With the estimate $\riem^{\f}=O(R^{-3})$, one gets a constant $C$ independent of the points and $Y$ such that $d_{\f}(x_1,x_2)\leq C\big(1+R(y)\big)^{-4-2\nu}\|Y\|_{C^{k,\alpha}_{\nu}}d_{\f}(x_1,x_2)$. 
 This can be made within two steps; first we call $\gamma_i$ the geodesic $t\mapsto\exp^{\f}_{x_i}\big(tY(x_i)\big)$:
  \begin{enumerate}
   \item control $d_{\f}(x_1,x_2)$ by $|p_{\gamma_1}(1)(Y(x_1))-p_{\gamma_2}(1)(Y(x_2))|^2 R(y)^{-3}$ (use \cite[Prop. 6.6]{bk});
   \item interpolate between the geodesics $\gamma_1$ and $\gamma_2$ via a minimizing geodesic $\alpha=\exp_{x_1}^{\f}(\cdot Z)$ joining $x_1$ to $x_2$ to get 
     \begin{equation*}
      |p_{\gamma_1}(1)(Y(x_1))-p_{\gamma_2}(1)(Y(x_2))|\leq Cd_{\f}(x_1,x_2)\|Y\|_{C^{k,\alpha}_{\nu}}(1+R(y)\big)^{-1-\nu}; 
     \end{equation*}
    the interpolation can be written for example as $\gamma_s(t)=\exp_{\alpha(s)}^{\f}\big[tY\big(\alpha(s)\big)\big]$. \cqfd 
  \end{enumerate}

 \subsubsection{The gauge and its consequences}

  We are now ready to put the metric $g_{\psi}$, pulled-back to $\C^2$ and extended smoothly, in so-called \textit{Bianchi-gauge} with respect to $\f$, at least outside a compact set. 
 We denote by $B^g=\delta^g+\tfrac{1}{2}d\tr^g$ the Bianchi operator of any Riemannian metric on $\C^2$. 
 Recall that $\Phi$ is a holomorphic identification between infinities of $\C^2/\Gamma$ and $X$, and that $\pi$ is the projection $\C^2\to \C^2/\Gamma$. 

 We assume that the parameter $\delta_1$ in Proposition \ref{prop_RFatinfty} is fixed $>\tfrac{1}{2}$. 
We state:
   \begin{prop} \label{prop_gauge}
  Let $(\alpha_2,\delta_2)\in(0,1)^2$ such that $\alpha_2<\alpha_1$, $\tfrac{1}{2}<\delta_2<\delta_1$, with $(\alpha_1,\delta_1)$ fixed in Proposition \ref{prop_RFatinfty}. 
  There exists a smooth diffeomorphism $\phi\in  \diff^{1,\alpha_2}_{\delta_2}$ such that
   \begin{equation}  \label{eq_gge}
    B^{\phi^*\tilde{\f}}\big((\Phi\circ\pi)^*g_{\psi}\big) = 0
   \end{equation}
  near infinity on $\C^2$, where $g_{\psi}$ stands for the metric of Proposition \ref{prop_RFatinfty}. 
  Moreover, $\phi$ can be chosen so as to commute with the action of $\Gamma$. 
 \end{prop}
 
Before proving Proposition \ref{prop_gauge}, we state the following consequence:

 \begin{crl}  \label{crl_rglrty}
  Let $\delta_2\in(\tfrac{1}{2},1)$ and $\phi$ as in Proposition \ref{prop_gauge}, and denote by $\tilde{\phi}$ the diffeomorphism induced on the infinity of $X$ by the diffeomorphism $\phi$ of Proposition \ref{prop_gauge}. 
  Then $(g_{\psi}-\tilde{\phi}^*\tilde{\f})\in C^{\infty}_{\delta_2+1}(X,\tilde{\phi}^*\tilde{\f})$.  
 \end{crl} 
  \prf \textit{of Proposition \ref{prop_gauge}}. 
 Recall that $\chi$ is the smooth nondecreasing cut-off function defined such that \eqref{eq_chi} holds. 
 Now, for $R_1$ big enough,
  \begin{equation*}
   g_{R_1}:= \chi(R-R_1)(\Phi\circ\pi)^*g_{\psi} + \big(1-\chi(R-R_1)\big)\f
  \end{equation*} 
 is a well defined metric on $X$, and is in the set $\met^{1,\alpha_1}_{\delta_1+1}(\f)$ of metrics $C^{1,\alpha_1}_{\delta_1+1}$ close to $\f$; 
 it is moreover $\Gamma$-invariant. 
 We want to find $\phi\in \diff^{1,\alpha_2}_{\delta_2}$ such that $B^{\phi^*\f}(g_{R_1}) = 0$, so that \eqref{eq_gge} will be verified on $\{R\geq R_1+1\}$. 
 Hence we are led to look at the operator 
  \begin{equation*}
   \begin{aligned} 
    \Psi\,:\, \diff^{2,\alpha_2}_{\delta_2}\times \met^{1,\alpha_2}_{\delta_2+1}(\f) \, &\longrightarrow \, C^{0,\alpha_2}_{\delta_2+1}(T^*X)   \\
                             (\phi,g)      \quad  \,\,\,       \qquad      \quad        &\longmapsto     \,           B^{\phi^*\f}(g_{R_1}),
   \end{aligned}
  \end{equation*}
where of course, $C^{0,\alpha_2}_{\delta_2+2}(T^*\C^2)$ denotes the space of $C^{0,\alpha_2}_{\delta_2+2}$ 1-forms on $\C^2$. 
Indeed, we want to use the implicit function theorem in a neighbourhood of $(\id,\f)$ to find $\phi$ \textit{once} $g$ is fixed, and we are interested in choosing $g=g_{R_1}$. 
Now, this latter $R_1$ can be arbitrarily big, provided that $\|g_{R_1}-\f\|_{C^{1,\alpha_2}_{\delta_2}}$ is arbitrarily small. 
But this is precisely the case, since $(g_{\psi}-\tilde{\f}) \in C^{1,\alpha_1}_{\delta_1}$ on $X$, and $\alpha_2<\alpha_1$, $\delta_2<\delta_1$ 
(whereas we would have no mean to claim that $\|g_{R_1}-\f\|_{C^{1,\alpha_1}_{\delta_1}}\to 0$ when $R_1$ goes to $\infty$). 

It only remains to observe that the derivative of $\Psi$ at $(\id,\f)$ with respect to its first argument is an isomorphism between 
$C^{2,\alpha_2}_{\delta_2}(T\C^2)$ and $C^{0,\alpha_2}_{\delta_2+2}(T^*\C^2)$. 
But this derivative is actually, after identifying $C^{0,\alpha_2}_{\delta_2+2}(T^*\C^2)$ to $C^{0,\alpha_2}_{\delta_2+2}(T\C^2)$, the rough Laplacian $(\nabla^{\f})^*\nabla^{\f}$ acting on vector fields, since $\f$ is Ricci-flat. 
It is precisely an isomorphism in our situation; until the end of this proof, since we do not use $\alpha_1$ and $\delta_1$ anymore, $\alpha$ and $\delta$ stand for $\alpha_2$ and $\delta_2$. 
 \begin{itemize}
  \item[$\bullet$] \textit{injectivity}: take $v\in C^{2,\alpha}_{\delta}$ so that $(\nabla^{\f})^*\nabla^{\f}v=0$. 
    For $t>0$, one can write the integration by parts
    \begin{equation*}  \label{eq_ipp}
     \int_{\{R\leq t\}} \langle v,(\nabla^{\f})^*\nabla^{\f}v\rangle_{\f} \vol^{\f} = \int_{\{R\leq t\}} \big|\nabla^{\f}v\big|^2\vol^{\f} -\int_{\{t=R\}} v\odot \nabla^{\f}v \vol^{\f|_{\{t=R\}}}.
    \end{equation*}
   The boundary term is a $O(t^{2-\delta-(\delta+1)})=O(t^{1-2\delta})$ (the spheres of radius $t$ have their volume in $t^2$ in Taub-NUT geometry), and hence $\nabla^{\f}v=0$ since $\delta>\tfrac{1}{2}$. 
   Such a $v$ is thus parallel, and since it tends to 0, we must have $v=0$. 
   Notice that since there is no critical value in $(0,1)$, the injectivity of $(\nabla^{\f})^*\nabla^{\f}$ still holds when $\delta\in (0,\tfrac{1}{2}]$. 
  \item[$\bullet$] \textit{surjectivity}: according to the theory of elliptic operators in weighted spaces, 
   the surjectivity of $(\nabla^{\f})^*\nabla^{\f}:C^{2,\alpha}_{\delta}\to C^{0,\alpha}_{\delta+2}$ amounts to the injectivity of $(\nabla^{\f})^*\nabla^{\f}$ on $C^{0,\alpha}_{1-\delta}$. 
   Thus if $w\in C^{0,\alpha}_{1-\delta}$ is in the kernel of $(\nabla^{\f})^*\nabla^{\f}$, we get by weighted elliptic estimates (see the techniques of \cite[App.]{bm}) that $w\in C^{2,\alpha}_{1-\delta}$. 
   Now we noticed that $(\nabla^{\f})^*\nabla^{\f}$ is injective on $C^{2,\alpha}_{1-\delta}$, and thus $w=0$. 
 \end{itemize}

Finally, we have to check that the vector field $Y$ such that $\Psi(\phi_Y,g_{R_1})=0$ given by the implicit function theorem is $\Gamma$-invariant; 
it is a direct consequence of its uniqueness (if chosen close enough to 0), and the $\Gamma$-invariance of the problem $\Psi(\phi_{\cdot},g_{R_1})=0$. 
The smoothness of $Y$, and hence that of $\phi_Y$, is purely local. 
\cqfd

\begin{rmk}
 We chose to work back on $\C^2$ with $\f$ instead of staying on $X$ with $\tilde{\f}$ to have exactly a rough Laplacian, and not $(\nabla^{\tilde{\f}})^*\nabla^{\tilde{\f}}+\ric^{\tilde{\f}}$. 
 Indeed, since $\ric^{\tilde{\f}}$ has no reason to be small where we extended $\f$ on $X$, we would not have been able to give an isomorphism statement for the latter operator. 
 It reduces nonetheless to a rough Laplacian near infinity; an alternate proof could thus consist in working with diffeomorphism of $X$ minus some compact set, and solve a Dirichlet problem as linearized problem. 
\end{rmk}

\prf \textit{of Corollary \ref{crl_rglrty}}. 
Take $\phi$ as in Proposition \ref{prop_gauge}. 
Then since $g_{\psi}$ is Ricci-flat at infinity, one has, outside of a compact subset of $X$,
 \begin{equation*}
  \left\{\begin{aligned}
          \ric^{g_{\psi}}   \quad                 &= 0,  \\
           B^{\tilde{\phi}^*\tilde{\f}}(g_{\psi}) &= 0.
         \end{aligned}
  \right.
 \end{equation*}
Set $\tilde{\F}:=\tilde{\phi}^*\tilde{\f}$ where it makes sense, and extend it as a Riemannian metric on $X$. 
Set also 
 \begin{equation*}
  \Phi^{\tilde{\F}}(h) = \ric(h)+(\delta^h)^* B^{\tilde{\F}}h
 \end{equation*}
for locally $C^2$ Riemannian metrics $h$, so that $\Phi^{\tilde{\F}}(\tilde{\F})=\Phi^{\tilde{\F}}(g_{\psi})=0$ near infinity. 
Set now $\vareps:=g_{\psi}-\tilde{\F}=(g_{\psi}-\tilde{\f})+(\tilde{\f}-\tilde{\F})\in C^{1,\alpha_2}_{1+\delta_2}$. 
Since $\Phi^{\tilde{\F}}$ is an operator of order 2, we can write schematically
 \begin{equation*}  
  0=\Phi^{\tilde{\F}}(\tilde{\F})-\Phi^{\tilde{\F}}(g_{\psi})=\big(d_{\tilde{\F}}\Phi^{\tilde{\F}}\big)(\vareps)+P\big(\vareps, \partial\vareps, \partial^2\vareps\big).
 \end{equation*}
where $P\big(\vareps, \partial\vareps, \partial^2\vareps\big)$ is some combination of $\vareps$, its first and second derivatives, which is at least quadratic, and with coefficients depending on $\F$. 
Now, in local coordinates, 
 \begin{equation}   \label{eq_PhiFh}
  \begin{aligned} 
   \big(\Phi^{\tilde{\F}}(&h)\big)_{ij} = -\frac{1}{2}\tilde{\F}^{kr}(\partial_{ij}h_{kr}-\partial_{ir}h_{jk}-\partial_{jk}h_{ir}+\partial_{kr}h_{ij})            \\
                         & +(\tilde{\F}^{kr}-h^{kr})\big(\partial_{ij}h_{kr}-\tfrac{1}{2}(\partial_{ir}h_{jk}-\partial_{jk}h_{ir})+\partial_{kr}h_{ij}\big)\\
                             & +\frac{1}{2}\big[(\partial_j h^{kr})(\partial_ih_{kr}+\partial_k h_{ir}-\partial_rh_{ik})                     
                                                -(\partial_i h^{kr})(\partial_jh_{kr}+\partial_k h_{jr}-\partial_rh_{jk})\big]                \\
                             & -\big(\Gamma^{\ell}_{jk}(h)\Gamma^k_{i\ell}(h)-\Gamma^{\ell}_{ij}(h)\Gamma^{k}_{k\ell}(h)\big)                 \\
                             &+S_{ij}\big[h_{j\ell}(\partial_ih^{\ell p})(\partial_kh_{rp}+\partial_rh_{pk}-\partial_ph_{rk})                 
                                          +h_{k\ell}(\partial_ih^{\ell p})(\partial_jh_{rp}+\partial_rh_{jp}-\partial_ph_{jr})\big]           \\
                             & -\big[\tilde{\F}^{kr}(\tfrac{1}{2}\partial_{\ell} h_{kr}-\partial_kh_{\ell r})
                                     +\tfrac{1}{2}h_{kr}\partial_{\ell} \tilde{\F}^{kr}+\tilde{\F}^{kr}h_{\ell m}\Gamma^{\ell}_{kr}(h)
                                     +\tilde{\F}^{kr}h_{km}\Gamma^m_{\ell r}(h)\big]
  \end{aligned}
 \end{equation}
where $S_{ij}$ means the symmetrization with respect to indexes $i$ and $j$ of its argument, and $\Gamma(h)$ the Christoffel symbols of $h$. 

The interest of this formula lies in the following: in $P\big(\vareps, \partial\vareps, \partial^2\vareps\big)$, 
 \begin{enumerate}
  \item the only occurrence of the second derivatives of $\vareps=g_{\psi}-\tilde{\F}$, which we denote by $\partial^2\vareps$ in \eqref{eq_PhiFh}, is via a  tensor of type $\vareps\odot \partial^2\vareps$, where $\odot$ is some algebraic operation with coefficients depending only on $g_{\psi}$ and $\tilde{\F}$;
  \item there is no appearance in \eqref{eq_PhiFh} of terms of type $\vareps\odot \vareps$, unless through a term of type $\vareps\odot \vareps\odot\partial\vareps$;
  \item the algebraic coefficients are controlled (for $\tilde{\F}$ say) in $C^{1,\alpha_2}$.
 \end{enumerate}
Let us sum up those three points by saying that 
 \begin{equation}  \label{eqn_lapleps}
  \big(d_{\tilde{\F}}\Phi^{\tilde{\F}}\big)(\vareps)+\vareps\odot \partial^2\vareps=\vareps\odot\vareps\odot P_1(\vareps,\partial\vareps)+
                                                                 \vareps\odot\partial\vareps\odot P_2(\vareps,\partial\vareps),
 \end{equation}
which provides us with control of $\big(d_F\Phi^{\F}\big)(\vareps)+\vareps\odot \partial^2\vareps$ in $C^{0,\alpha_2}_{2\delta_2+3}$. 
On the other hand, since $\tilde{\F}$ is Ricci-flat in a neighbourhood of infinity, the derivative $d_{\tilde{\F}}\Phi^{\tilde{\F}}$ there is well-known: it is nothing but $\tfrac{1}{2}\mathscr{L}_{\tilde{\F}}$, where $\mathscr{L}_{\tilde{\F}}$ denotes the Lichnerowicz Laplacian associated to $\tilde{\F}$. 
Now since $\vareps\in C^{1,\alpha_2}_{1+\delta_2}$, the perturbed Laplacian
 \begin{equation*}
  \EuScript{L}_{\vareps}:\hat{\vareps} \longmapsto \frac{1}{2}\mathscr{L}_{\tilde{\F}}\hat{\vareps}+\vareps\odot\partial^2\hat{\vareps}
 \end{equation*}
(where we emphasize that the bareheaded $\vareps$ is seen as \textit{fixed}) is $C^{1,\alpha_2}_{1+\delta_2}$ close to $\tfrac{1}{2}\mathscr{L}_{\tilde{\F}}$, and in particular is elliptic near infinity. 
The weighted estimates are thus the same for this operator as for $\tfrac{1}{2}\mathscr{L}_{\tilde{\F}}$ (at least near infinity). 
Since $\vareps\in C^{1,\alpha_2}_{1+\delta_2}$, we get from \eqref{eqn_lapleps}, which we rewrite as $\EuScript{L}_{\vareps}\vareps\in C^{0,\alpha_2}_{2\delta_2+3}$, 
that $\vareps\in C^{2,\alpha_2}_{1+\delta_2}$, near infinity and hence on $X$, since it is smooth. 

Plugging back this estimate into \eqref{eqn_lapleps}, we get $\EuScript{L}_{\vareps}\vareps\in C^{1,\alpha_2}_{2\delta_2+3}$, with $\EuScript{L}_{\vareps}$ $C^{2,\alpha_2}_{\delta_2+1}$ close to $\tfrac{1}{2}\mathscr{L}_{\tilde{\F}}$. 
It follows that $\vareps\in C^{3,\alpha_2}_{1+\delta_2}$ near infinity, and thus on the whole $X$. 
An immediate inductive repetition of these arguments provides finally that $\vareps\in C^{\infty}_{\delta_2+1}$.
\cqfd

~

\paragraph{Reformulation of Corollary \ref{crl_rglrty}.} 
Set now (on $\C^2$) $\overline{y}_j=\phi^*(y_j)$ and $\overline{\eta}=\phi^*(\eta)$, and so on, so that near infinity on $X$, $\tilde{\F}= \overline{V}(d\overline{y}_1^2+d\overline{y}_2^2+d\overline{y}_3^2)+\overline{V}^{-1}\overline{\eta}^2$. 
Forgetting the projection $\pi$ and the identification $\Phi$, since $I_1$ is parallel for $g_{\psi}$, in view of Corollary \ref{crl_rglrty}, it is a rather simple exercise to write
 \begin{equation*}
  g_{\psi}= \overline{V}(d\overline{y}_1^2+d\overline{y}_2^2+d\overline{y}_3^2)+\overline{V}^{-1}\overline{\eta}^2+\vareps
 \end{equation*}
with $\vareps\in C^{\infty}_{\delta_2+1}$, and $(I_1\overline{V}d\overline{y}_1-\overline{\eta})$ and $(I_1d\overline{y}_2-d\overline{y}_3)\in C^{\infty}_{\delta_2+1}(\tilde{\F})$.

\subsection{Making the metric Ricci-flat on the whole $X$.}  \label{sctn_ccl_pgrm}

We perform now the very last step of our program. 
For this we state a general Calabi-Yau theorem for Kähler ALF 4-manifolds of dihedral type --terminology is made clear right after the statement--, which we prove in next part:
 \begin{thm} \label{thm_CY_ALF}
  Let $(Y,g_Y,J^Y,\omega_Y)$ be an ALF Kähler 4-manifold of dihedral type. 
  Given a $C^{\infty}_{loc}$ function $f$ such that there exists some positive $\beta\in(0,1)$ such that $(\nabla^{g_Y})^{\ell} f=O(R^{-2-\ell-\beta})$ for all $\ell\geq 0$, 
  there exists a unique $\varphi\in C^{\infty}_{loc}$ such that $(\nabla^{g_Y})^{\ell} \varphi=O(R^{-\ell-\beta})$ for all $\ell\geq 0$ and 
  \begin{equation} \label{MAeqn}
   \big(\omega_Y+i\ddbar\varphi\big)^2=e^f\omega_Y^2.
  \end{equation}
 \end{thm}

~

\paragraph{4-manifolds of ALF dihedral type.} 
Before applying this theorem to our construction, we shall define the terms it involves. 
In particular, let us make precise the statement "$(Y,g_Y,J^Y,\omega_Y)$ is an ALF Kähler 4-manifold of dihedral type". 
We say so if $(Y,g_Y,J^Y,\omega_Y)$ is a complete Kähler manifold of real dimension 4 such that:
 \begin{itemize}
  \item[$\bullet$] outside a compact subset, $Y$ admits a holomorphic two-sheeted covering space $(\EuScript{Y},J^{\EuScript{Y}})$ which is the total space of an $S^1$-fibration $\varpi=(y_1,y_2,y_3)$ over $\R^3$ minus some ball; 
  \item[$\bullet$] the fibration is covariant under the $\Z_2$-action on $\EuScript{Y}$: if $\alpha$ generates the $\Z_2$-action, then $\varpi\big(\tau(y)\big)=-\varpi(y)$ for any $y\in\EuScript{Y}$;
  \item[$\bullet$] $g_Y$ is asymptotic to some metric $h$ adapted to the fibration; more precisely, there exists $L>0$ and an $S^1$-invariant 1-form $\mu$ such that 
   if $\tfrac{L}{2\pi}T$ is the infinitesimal generator of the $S^1$-action on the fibers, then $\mu(T)=1$, and 
    \begin{equation*}
     h :=\varpi^*g_{\R^{3}}+\mu^2 \quad \text{on} \quad\EuScript{Y},
    \end{equation*}
   then for some $\nu>0$, if we still denote by $h$ its push-forward to $Y$,
    \begin{equation*}
     (\nabla^{h})^\ell(g-h)=O(\rho^{-\nu-\ell}) \quad \text{for all } \ell\geq0.
    \end{equation*}
   We also assume the following decay condition on the connection form $\mu$:
    \begin{equation*} 
     D^{\ell}\mu=O(\rho^{-\nu-1-\ell}) \quad\text{for all } \ell\geq1, 
    \end{equation*} 
   where $D$ is the standard connection on $\R^{3}$ and $\rho$ is the radius on $\R^3$, transposed and extended smoothly on $Y$. 
  \item[$\bullet$] the complex structure $J^Y$ and $\varpi$ are compatible in the sense that: $D^{\ell}\big(J^{\EuScript{Y}}dx_1-\mu\big)=O(R^{-\ell-\tau})$ and $D^{\ell}\big(J^{\EuScript{Y}}dx_2-dx_3\big)=O(R^{-\ell-\tau})$ for all $\ell\geq0$. 
 \end{itemize}
One can check that if $(Y,g_Y,J^Y,\omega_Y)$ verifies this definition, then $(\nabla^{g_Y})^{\ell}\riem^{g_Y}=O(\rho^{-\nu-2-\ell})$ for all $\ell\geq0$. 

~

\paragraph{Application of Theorem \ref{thm_CY_ALF}.}
Set $Y=X$, $J^Y=I^X_1$, $g_Y=g_{\psi}$, and take the fibration induced by $\phi^*(y_1,y_2,y_3)$ on $(\C^2\backslash B)/<\zeta_k>$, with $\phi$ the diffeomorphism involved in the gauge process (section \ref{sctn_gge}). 
This provides $(X,g_{\psi},I^X_1,\omega_{\psi})$ the structure of an ALF 4-manifold of dihedral type. 
The parameter $L$, the geometric interpretation of which is the length of the fibers at infinity, becomes $\tfrac{\pi\sqrt{2/m}}{k}$. 
Thus set $f = \log\big(\tfrac{\Omega_g}{\vol^{g_{\psi}}}\big)$, so that $\Omega_g= e^{f} \vol^{g_{\psi}}$; since $f$ is smooth and has compact support, it verifies the hypotheses of Theorem \ref{thm_CY_ALF} for any $\beta\in (0,1)$. 

The theorem now gives us a $\varphi$ which is in $ C^{\infty}_{\beta}$ for any $\beta\in(0,1)$ and such that 
 \begin{equation*}
  \big(\omega_{\psi}+i\ddbar\varphi\big)^2= e^f\omega_{\psi}^2 = 2\Omega_g,
 \end{equation*}
and $\omega_{RF}:=\omega_{\psi}+i\ddbar\varphi$ is Kähler by the usual arguments (equivalent to $\omega_{\psi}$ at infinity, nowhere-vanishing determinant). 
The associated Riemannian metric, $g_{RF}$ say, is thus Kähler for $I^X_1$ and Ricci-flat, and $(g_{RF}-g_{\psi})\in C^{\infty}_{\beta+2}$ for any $\beta\in(0,1)$. 

~

\paragraph{Conclusion of the program.} 
The Ricci-flat metric $g_{RF}$ is also hyperkähler; indeed, set $\omega_j= g(I^X_j\cdot,\cdot)$, $j=2,3$, and denote by $\theta$ the holomorphic (2,0)-form $\omega_2+i\omega_3$. 
Then $\theta\wedge\overline{\theta}= 4 \Omega_g$. 
Hence $\theta$ has constant norm with respect to $g_{RF}$ and is thus parallel ; so are $\omega_2$ and $\omega_3$. 
Define endomorphisms $J^X_j$ of $TX$ by $\omega_j= g_{RF}(J^X_j\cdot,\cdot)$, $j=2,3$. 
Elementary manipulations tell us these are almost-complex structures, that $g_{RF}$ is hermitian for them, and that $J^X_1J^X_2J^X_3=-1$ if one sets $J^X_1=I^X_1$. 
They are automatically parallel, hence integrable: $(g_{RF},J^X_j)$ is hyperkähler. 

To sum our construction up, we say for example that $g_{RF}$ is Kähler with respect to $I_1^X$, hyperkähler with volume form $\Omega_g$, and $(g_{RF}-\f)\in C^0_{\delta_1+2}\cap C^1_{\delta_1+1}$. 
Now $g_{RF}$ may depend on our intermediate steps, like Propositions \ref{prop_omegam} and \ref{prop_RFatinfty} ; suppose then $g_{RF}'$ is obtained with the same program, but with different parameters ; in particular, call $\delta_1'$ the relevant one. 
Following the program nevertheless, we see that $\omega_{RF}-\omega'_{RF}$ can be written as $dd^c \nu$ for some $\nu\in C^2_{\min(\delta_1,\delta_1')}$. 
Moreover $g_{RF}$ and $g_{RF}'$ have the same volume form; as a consequence, $\Delta_{g_{RF}}\nu$ has constant sign, and thus $\nu\equiv0$, i.e. $g_{RF}=g'_{RF}$

A consequence of this uniqueness property is that $(g_{RF}-\tilde{\f})\in C^0_{\delta+2}\cap C^1_{\delta+1}$ \textit{for any} $\delta\in(0,1)$; the naive prototype $(\Phi\circ\pi)_*\f$ on $X$ is thus a rather good approximation to the desired Ricci-flat metric. 
Despite this though, let us observe that we have to make some detour by a diffeomorphic perturbation of $\f$ to perform the last step of our construction. 
This concludes the program we fixed at the beginning of this part, and completes the proof of Theorem \ref{thm_gRF}. 

~

\section{Proof of Theorem \ref{thm_CY_ALF}}  \label{part_prf_CY_ALF}

We shall use for proving Theorem \ref{thm_CY_ALF} a classical method in the study of Monge-Ampère equations: the \textit{continuity method}. 
This method was suggested by E. Calabi for the solution of the complex Monge-Ampère equation on compact Kähler manifolds. 
Since the successful use by Yau \cite{yau} of this method, it has been adapted to different non-compact settings; 
let us quote here the version by Joyce \cite[ch. 8]{joy} for ALE manifolds, which greatly inspired ours. 
We also refer the reader to Tian and Yau's seminal work \cite{tian-yau1, tian-yau2}, which pioneered the research on generalizing Calabi-Yau theorem to non-compact manifolds. 
We shall mention a result by Hein \cite[Prop. 4.1]{hein} too, very similar to ours if taking 
Hein's parameter $\beta$, but dealing with less precise asymptotics. 

We hence start this part by describing the method, and follow by the analytic work (in particular, a priori estimates) it requires. 

\subsection{The continuity method.}  \label{sect_cntnty_mthd}

The principle of this method is not to solve directly the Monge-Ampère equation (\ref{MAeqn}), but to solve a one-parameter family of such equations in which the right-hand side term moves from $\omega_Y^2$ to $e^f\omega_Y^2$. 
Concretely, we consider the equations 
 \begin{equation}  \tag{$E_t$}
  \big(\omega_Y+i\ddbar\varphi_t\big)^m=e^{tf}\omega_Y^m
 \end{equation}
for $t\in[0,1]$. 
Now set $S:=\{t\in[0,1]|\, (E_t)\text{ has a unique solution } \varphi_t\in C^{\infty}_{\beta}\}$; 
the weighted space involved in this definition is defined on $Y$ in total analogy with those of section \ref{sctn_RFnrinfty} (replace $g_m$ by $g_Y$, $\tilde{R}$ by $\rho$, and so on), 
and we could express Theorem \ref{thm_CY_ALF} by writing $f\in C^{\infty}_{\beta+2}$ and $\varphi\in C^{\infty}_{\beta}$. 
  
Up to the uniqueness of the solution $0$ of $(E_0)$, it is obvious that $0\in S$. 
We shall then prove: 
 \begin{enumerate}
  \item the set $S$ is closed (section \ref{closedness}); 
  \item the set $S$ is open (section \ref{sect_openness}). 
 \end{enumerate}
Theorem \ref{thm_CY_ALF} then follows from an immediate connectedness argument.

In order to prove the openness of $S$, one observes that the linearisation of the Monge-Ampère operator is (nearly) a Laplacian, close enough to $\Delta_{g_Y}$; this allows one to use the isomorphisms induced between some Hölder spaces, plus ellipticity of such an operator. 
This is done in section \ref{sect_openness}. 

On the other hand, in order to prove the closedness of $S$ (the hard part), we need some compactness properties for solutions of a family of $(E_t)$, and this we get by establishing a priori estimates on such solutions. 
This is done in paragraph \ref{prgph_clsdnss}.

The  easiest part is the uniqueness of the solution of any $(E_t)$, and we shall deal with it now as a warm-up.
 \begin{lem}
  Let $t\in[0,1]$, and let $\varphi_1$, $\varphi_2$ be two solutions of $(E_t)$. 
  Assume that $\varphi_1$ and $\varphi_2\in C^{\infty}_{\beta}$. 
  Then $\varphi_1=\varphi_2$.
 \end{lem}
\prf. This is essentially the same argument as that used to prove the uniqueness of $g_{RF}$ above, namely if for some $t$, $\varphi_1$ and $\varphi_2$ are $C^{2}_{\beta}$ solutions of $(E_t)$, then $\omega_{\varphi_1}=\omega_Y+i\ddbar\varphi_1$ is a Kähler form equivalent to $\omega_Y$, 
and denoting by $\Delta_1$ the Laplacian of its associated metric, one gets that $\Delta_1(\varphi_1-\varphi_2)$ has constant sign. 
Since $\varphi_1-\varphi_2$ tends to 0 at infinity, we thus have $\varphi_1=\varphi_2$. 
\cqfd

\subsection{Closedness of $S$: a priori estimates.} \label{closedness}

 \subsubsection{$C^0$ estimates.}  \label{c0_estimate}

We are proving in this paragraph an a priori estimate on the $C^0$ norm of a solution of a Monge-Ampère equation on $Y$. 

The estimates can already be deduced from \cite[Prop. 4.1]{hein}; 
we give our own proof though, because \textit{weighted $C^0$ estimates} will be established below (\S \ref{c0delta_estimate}) in a very similar but slightly more complicated way. 

The techniques used here are quite close to the ones used by Yau in the compact case, namely a recursive use of integration by parts giving a variation of Moser's iteration adapted to the Monge-Ampère equation, 
and their adaptation by Joyce for his version of the Calabi-Yau theorem on ALE manifolds. 
Nonetheless, because ALF geometry is quite different from ALE geometry, we write in detail the manipulations that need to be adapted to our framework. 
In particular, because of different measures on both sides of a Sobolev inequality (Lemma \ref{lem_sob_inj}) we make a crucial use of, our system of weights in the following integrals is not the same as Joyce's. 
 
We start with the following :
 \begin{lem} \label{lem_ibp_formula}
  Let $f\in C^0_{\beta+2}$, $\beta\in(0,1)$, and $\varphi\in C^3_{\gamma}$, $0<\gamma\leq\beta$, such that $\big(\omega_Y+i\ddbar\varphi\big)^2=e^{f}\omega_Y^2$. Then for any $p>\gamma^{-1}$, $p>2$, 
  \begin{equation*}
   \int_Y \big|\partial|\varphi|^{p/2}\big|^2\omega^2_Y\leq \frac{p^2}{2(p-1)}\int_Y |\varphi|^{p-2}\varphi (e^f-1)\omega^2_Y,
  \end{equation*}
 where both integrals are finite.
 \end{lem}

Notice that in this lemma $f$ plays the role of $tf$ in our problem, but this is not an issue since a $C^0_{\beta+2}$ bound on $f$ gives uniform $C^0_{\beta+2}$ bounds on the $tf$.

~

\noindent\prf \textit{of Lemma \ref{lem_ibp_formula}}. 
First, we set $T=\omega_Y+\omega_{\varphi}$, with $\omega_{\varphi}=\omega_Y+i\ddbar\varphi$. 
Picking $p >\gamma^{-1}$ and $p>2$, we claim that 
 \begin{equation*}
  \int_{Y} d\big(\varphi|\varphi|^{p-2} d^c\varphi \wedge T\big)=0.
 \end{equation*}
Indeed, for any real number $R$ big enough, if $B_R$ denotes the set $\{\rho<R\}$ and $S_R$ its boundary, Stokes' theorem asserts that
 \begin{equation*}
  \int_{B_R} d\big(\varphi|\varphi|^{p-2} d^c\varphi \wedge T\big)=\int_{S_R} \varphi|\varphi|^{p-2} d^c\varphi \wedge T.
 \end{equation*}
Now since the volume of $S_R$ according to $g$ is $O(R^2)$, since $T$ is bounded with respect to $g_Y$ and since $\varphi|\varphi|^{p-2} d^c\varphi=O(R^{-p\gamma-1})$ on $S_R$, we get that the right-hand-side term is $O(R^{1-p\gamma})$, and hence goes to 0 as $R$ goes to $\infty$ by our choice of $p>\gamma^{-1}$. 
A little computation yields
 \begin{equation*}
  \frac{1}{2} d\big(\varphi|\varphi|^{p-2} d^c\varphi \wedge T\big)=\varphi|\varphi|^{p-2} i\ddbar\varphi \wedge T
                                                                    +(p-1)|\varphi|^{p-2} i\partial\varphi\wedge\dbar\varphi\wedge T
 \end{equation*}
because $T$ is closed. 
But on the one hand, $i\ddbar\varphi \wedge T=(\omega_{\varphi}-\omega_Y)\wedge T=(e^f-1)\omega_Y^2$ and on the other hand $|\varphi|^{p-2} i\partial\varphi\wedge\dbar\varphi=\tfrac{4}{p^2}i\partial\big(|\varphi|^{p/2}\big)\wedge\dbar\big(|\varphi|^{p/2}\big)$, hence
 \begin{equation*}
  \int_Y i\partial\big(|\varphi|^{p/2}\big)\wedge\dbar\big(|\varphi|^{p/2}\big)\wedge T=\frac{p^2}{2(p-1)}\int_Y \varphi|\varphi|^{p-2}(1-e^f)\omega_Y^2.
 \end{equation*}
Finally, observe that $\tfrac{i\alpha\wedge\overline{\alpha}\wedge\omega_{\varphi}}{\omega_{Y}^2}=\tfrac{1}{2}e^f|\alpha|_{\omega_{\varphi}}^2\geq 0$, and that $i\alpha\wedge\overline{\alpha}\wedge\omega_Y=\tfrac{1}{2}|\alpha|_{g_Y}^2\omega_Y^2$ for any (1,0)-form $\alpha$, and conclude by setting $\alpha=\partial\big(|\varphi|^{p/2}\big)$. 
The finiteness of the integrals in play merely comes from our choice of $p$. 
\cqfd

~

Next, we want to use the inequality we have just proved with some Sobolev embedding to estimate recursively $\|\varphi\|_{L^{p\vareps}_{d\lambda}}$ in terms of $\|\varphi\|_{L^{p}_{d\lambda}}$, with some fixed $\vareps>1$, and some positive measure $d\lambda$. 
The Sobolev embedding states as : 
 \begin{lem} \label{lem_sob_inj}
  There exists a constant $C_{S}$ such that for any $u\in H^1_{loc}$ (hence in $L^4_{loc}$) such that the $\int_{Y} |u|^{4}\rho^{-1}\vol^{g_Y}$ is finite, one has
   \begin{equation} \label{eqn_sob_inj}
   \Big(\int_{Y} |u|^{4}\rho^{-1}\vol^{g_Y}\Big)^{1/4} \leq C_{S}^{1/2} \Big(\int_Y  |du|_g^2 \vol^{g_Y} \Big)^{1/2}. 
   \end{equation}
 \end{lem}
This inequality can be compared to \cite[Prop. 3.2]{hein}; we postpone its proof at the end of the present paragraph. 

First we initiate the induction by estimating $\|\varphi\|_{L^{p_0}_{d\lambda}}$ (or $\|\varphi\|_{L^{p_0\vareps}_{d\lambda}}$) for some $p_0$ independent of $\varphi$. 
It will be relevant to take $\vareps=2$, and $d\lambda=\rho^{-1}\vol^{g_Y}$ in the following. 
Our initial estimation states:
 \begin{lem}  \label{lem_initiate}
  Fix some $p_0>2$, $p_0>\gamma^{-1}$. 
Then under the assumptions of Lemma \ref{lem_ibp_formula} there exists $C$ only depending on $\beta$, $\|f\|_{C^0_{\beta+2}}$, $p_0$ and $g$ such that $\|\varphi\|_{L^{p_0\vareps}_{d\lambda}}\leq C$.
 \end{lem}
\prf. Apply inequality \eqref{eqn_sob_inj} to $u=|\varphi|^{p_0/2}$ to get:
 \begin{equation*}
   \Big(\int_{Y} |\varphi|^{2p_0}\rho^{-1}\vol^g\Big)^{1/2} 
                                                           \leq C_{S} \int_Y  \big|d|\varphi|^{p_0/2}\big|_g^2 \vol^g . 
 \end{equation*}
By Lemma \ref{lem_ibp_formula}, we have
 \begin{equation*}
  \Big(\int_{Y} |\varphi|^{2p_0}\rho^{-1}\vol^{g_Y}\Big)^{1/2} 
                                                           \leq \frac{p_0^2C_{S}}{2(p_0-1)}\int_Y |\varphi|^{p_0-1} \big|e^f-1\big|\vol^{g_Y}.
 \end{equation*}
Write $\big|e^f-1\big|=\big|e^f-1\big|^{a+b}$ with $a=\tfrac{p_0-1}{2(\beta+2)p_0}$ and $b=1-a$, and apply Hölder inequality to the right-hand-side term of the inequality above with exponents $s=\tfrac{2p_0}{p_0-1}$ and $t=\big(1-\tfrac{1}{s}\big)^{-1}=\tfrac{2p_0}{1+p_0}$. 
This yields:
 \begin{align*}
  \Big(\int_{Y} |\varphi|^{2p_0}\rho^{-1} \vol^{g_Y} \Big)^{1/2} 
      &\leq \frac{p_0^2C_{S}}{2(p_0-1)} \Big(\int_Y|\varphi|^{s(p_0-1)}\big|e^f-1\big|^{as}\Big)^{1/s} \Big(\int_Y\big|e^f-1\big|^{bt}\Big)^{1/t} \\
     =& \frac{p_0^2C_{S}}{2(p_0-1)}\Big(\int_Y|\varphi|^{2p_0}\big|e^f-1\big|^{\tfrac{1}{(\beta+2)}}\Big)^{1/s}  
                                     \Big(\int_Y\big|e^f-1\big|^{bt}\Big)^{1/t}.
 \end{align*}
Noticing that $|e^f-1|=O(\rho^{-(\beta+2)})$ (and more precisely that at any point, $|e^f-1|\leq e^{\|f\|_{C^0}} \|f\|_{C^0_{\beta+2}} \rho^{-(\beta+2)}$), we get that $\int_Y|\varphi|^{2p_0}|e^f-1|^{\tfrac{1}{(\beta+2)}}\leq C\int_{Y} |\varphi|^{2p_0}\rho^{-1}$ for some $C$ depending only on the parameters announced. 
Moreover, we have that $(\beta+2)bt=\tfrac{2p_0}{1+p_0}\big((\beta+2)-\tfrac{p_0-1}{2p_0}\big)>3$ because $p_0>\gamma^{-1}\geq\beta^{-1}$, so that $\int_X\big|e^f-1\big|^{bt}$ is finite and equal to some constant, $K^t$ say, also independent of $\varphi$. 
So far, we obtain that:
 \begin{equation*}
  \|\varphi\|_{L^{2p_0}_{d\lambda}}^{p_0}\leq  \frac{p_0^2CC_{S}}{2(p_0-1)} K \|\varphi\|_{L^{2p_0}_{d\lambda}}^{p_0-1},
 \end{equation*}
and the conclusion follows. 
\cqfd

~

We fix now $p_0=\tfrac{2}{\gamma}$. 
Using the same techniques, we can prove a recursive control on $\|\varphi\|_{_{L^{2p}_{d\lambda}}}$ from $\|\varphi\|_{_{L^{p}_{d\lambda}}}$:
 \begin{lem}  \label{lem_propagate}
  Under the assumptions of Lemma \ref{lem_ibp_formula}, there exists a constant $C_1$ only depending on $\beta$, $\gamma$, $\|f\|_{C^0_{\beta+2}}$ and $g_Y$ such that 
  for any $p\geq p_0$, $\|\varphi\|_{L^{2p}_{d\lambda}}^{p}\leq C_1 p\|\varphi\|_{L^{p}_{d\lambda}}^{p-1}$.
 \end{lem}
\prf. The idea is the same as for the previous lemma: apply inequality \eqref{eqn_sob_inj} to $u=|\varphi|^{p/2}$, Lemma \ref{lem_ibp_formula} and Hölder inequality (with well-chosen exponents). 
The first two give, for $p>\gamma^{-1}$, $p>2$:
 \begin{equation*}
  \Big(\int_Y |\varphi|^{2p}\rho^{-1}\vol^{g_Y}\Big)^{1/2} 
                                                           \leq \frac{p^2C_{S}}{2(p-1)}\int_Y |\varphi|^{p-1} \big|e^f-1\big|\vol^{g_Y}.
 \end{equation*}
Apply now Hölder inequality with exponents $s=\tfrac{p}{p-1}$, $t=p$ and weights $a=\tfrac{(p-1)}{(\beta+2)p}=\tfrac{1}{(\beta+2)s}$ and $b=1-a$ to write:
 \begin{align*}
  \Big(\int_Y |\varphi|^{2p}\rho^{-1} & \vol^{g_Y}\Big)^{1/2}\\
                           \leq & \frac{p^2C_{S}}{2(p-1)}\Big(\int_Y |\varphi|^{p} \big|e^f-1\big|^{\tfrac{1}{(\beta+2)}}\vol^{g_Y}\Big)^{(p-1)/p}
                                                          \Big(\int_Y\big|e^f-1\big|^{bp}\vol^{g_Y}\Big)^{1/p}.
 \end{align*}
But $|e^f-1|^{\tfrac{1}{(\beta+2)}}\leq C\rho^{-1}$ for some constant depending only on $f$. 
Moreover, $(\beta+2)pb=p(\beta+1)+1$, which is $>3$ as soon as $p>\tfrac{2}{\beta+1}$; this is actually automatic since $p>\gamma^{-1}\geq\beta^{-1}$ and $\beta\in(0,1)$. 
Since furthermore $b=b(p)$ tends to $\tfrac{\beta+1}{\beta+2}:=b_{\infty}>0$ when $p$ goes to $\infty$, $\big(\int_Y |e^f-1 |^{bp}\vol^g\big)^{1/p}$ tends to $\sup_{Y}|e^f-1|^{b_{\infty}}$ when $p$ goes to $\infty$, so we can claim: for all $p\geq p_0$, $\big(\int_Y |e^f-1 |^{bp}\vol^{g_Y}\big)^{1/p}\leq K$ for some $K$ depending only on $f$ (and which we could evaluate in terms of $\|f\|_{C^0_{\beta+2}}$ only). 
Finally, we get that for all $p\geq p_0$:
 \begin{align*}
  \Big(\int_Y |\varphi|^{2p}\rho^{-1}\vol^{g_Y}\Big)^{1/2} 
                           &\leq \frac{p^2C_{S}C^{(p-1)/p}K}{2(p-1)}\Big(\int_Y |\varphi|^{p}\rho^{-1} \vol^{g_Y}\Big)^{(p-1)/p}\\
                           &\leq C_1 p \Big(\int_Y |\varphi|^{p}\rho^{-1} \vol^{g_Y}\Big)^{(p-1)/p}
 \end{align*}
with $C_1=(1+C)KC_S$ say, which only depends on the parameters announced. 
\cqfd

~

We fix now $p_0=2\gamma^{-1}$ in Lemma \ref{lem_initiate}.
 \begin{lem}  \label{lem_induction}
  Under the assumptions of Lemma \ref{lem_ibp_formula}, there exist two constants $Q_0$ and $C_2$ depending only on $\beta$, $\gamma$, $\|f\|_{C^0_{\beta+2}}$ and $g_Y$ such that 
  for any $q\geq 2p_0$, $\|\varphi\|_{L^{q}_{d\lambda}}\leq Q_0(C_2 q)^{-2/q}$.
 \end{lem}
Letting $q$ go to $\infty$, we get the $C^0$ a priori estimate for $\varphi$ we are seeking :
 \begin{prop}  \label{prop_c0_estimate}
  Under the assumptions of Lemma \ref{lem_ibp_formula}, there exists a constant $Q_0=Q_0(\beta,\gamma,\|f\|_{C^0_{\beta+2}},g_Y)$ such that $\|\varphi\|_{C^0}\leq Q_0$.
 \end{prop}

~

We conclude this paragraph by the proof of Lemma \ref{lem_sob_inj}:

~

\prf \textit{of Lemma \ref{lem_sob_inj}}: 
We shall prove that there exist two constants $C_1$ and $C_2$ such that for any $u$ as in the statement of the lemma, 
 \begin{equation}  \label{eqn_sob_ineq}
  \Big(\int_{Y} |u|^{4}\rho^{-1}\vol^{g_Y}\Big)^{1/4} \leq C_1 \Big(\int_Y  |du|_{g_Y}^2 \vol^{g_Y} + \int_Y |u|^2 \rho^{-2}\vol^{g_Y}\Big)^{1/2}
 \end{equation}
and 
 \begin{equation*}
  \int_Y |u|^2 \rho^{-2}\vol^{g_Y} \leq C_2 \int_Y  |du|_{g_Y}^2 \vol^{g_Y};
 \end{equation*}
the Lemma will then follow at once. 
Let us start by the latter inequality. 
We prove it first for $u$ with have compact support in $\{\rho\geq\rho_0\}$, where $\rho_0$ is chosen so that we have the two-sheeted cover $\EuScript{Y}$ described in section \ref{sctn_ccl_pgrm} over $\{\rho\geq\rho_0\}$; 
this way, we can replace $Y$ by $\EuScript{Y}$, $g_Y$ by $h$ and $u$ by its pull-back. 
Decompose $u$ as $u_{\perp}+u_{0}$, with $u_{0}(x)$ the mean value of $u$ along the fiber of $\varpi$ passing by $x$. 
This makes $u_0$ a compactly supported function on $\R^3\backslash B$ ($B$ the unit ball). 
Take spherical coordinates $(\rho, \theta,\phi)$ on $\R^3$. 
Then 
 \begin{equation*}
  0=\int_{\rho\geq1}\partial_{\rho}(u_0^2\rho) d\rho \vol^{S^2}=2\int_{\rho\geq1} u_0\partial_{\rho}(u_0)\rho d\rho \vol^{S^2}
                                                                                    +\int_{\rho\geq1} u_0^2 d\rho \vol^{S^2}
 \end{equation*}
which we rewrite as
 \begin{align*}
  \int_{\rho\geq1} u_0^2 \rho^{-2}\vol^{\R^3} &= -2\int_{\rho\geq1} u_0\partial_{\rho}(u_0)\rho^{-1}\vol^{\R^3}         \\
                                                  &\leq 2\Big(\int_{\rho\geq1} u_0^2 \rho^{-2}\vol^{\R^3}\Big)^{1/2}
                                                          \Big(\int_{\rho\geq1} \partial_{\rho}(u_0)^2\vol^{\R^3}\Big)^{1/2}
 \end{align*}
(Cauchy-Schwarz inequality), i.e. $\int_{\rho=1}^{\infty} u_0^2 \rho^{-2}\vol^{\R^3}\leq 4\int_{\rho=1}^{\infty} |du_0|^2\vol^{\R^3}$. 

Now for the component $u_{\perp}$, after pulling it back on $\EuScript{Y}$, one has $\int_{\EuScript{Y}} u_{\perp}^2\rho^{-2}\vol^{h}=\int_{\rho\geq\rho_0}\rho^{-2}\vol^{\R^3}\int_{\text{fiber}}u_{\perp}^2 \mu$.
Since $u_{\perp}$ has zero mean along the fibers, $\int_{\text{fiber}}u_{\perp}^2 \mu\leq c \int_{\text{fiber}}d|u_{\perp}|^2_{h} \mu$ for some $c$ independent of the fiber, and thus 
$\int_{\EuScript{Y}} u_{\perp}^2\rho^{-2}\vol^{h}\leq c\int_{\EuScript{Y}} |du_{\perp}|_h^2\rho^{-2}\vol^{h}$; we even have a $\rho^{-2}$ factor in the RHS, which we can loosely get rid of. 
Adding these estimations on separate components, we get that $\int_{\EuScript{Y}} u^2\rho^{-2}\vol^{h}\leq C\int_{\EuScript{Y}} |du|_h^2\vol^{h}$, which easily transposes on $Y$ with $g_Y$. 
Now, to extend this Hardy type inequality to a general $u$ as in the statement, one may proceed as in \cite[\S 1.4.1]{auv1}; the only point to be noticed is that the only integrable constant on $Y$ for $\rho^{-1}\vol^{g_Y}$ is 0 (this replaces a zero mean assumption).

We now prove the Sobolev type inequality \eqref{eqn_sob_ineq}. 
Since such an inequality is known on compact manifolds, we can assume that $u$ has compact support on $\{\rho\geq\rho_0\}$. 
We pull-back $u$ to $\EuScript{Y}$ once more, and split it again into $u_0+u_{\perp}$. 
We easily get the inequality on $u_{\perp}$ (and even a much better one) by using the standard Sobolev embedding on the circle. 
Now for the component $u_{\perp}$, we write $\R^3\backslash B=\bigcup_{\ell\geq0} \mathcal{A}_{\ell}$, where $\mathcal{A}_{\ell}$ is the annulus $\{2^{\ell}\leq \rho\leq 2^{\ell+1}\}$. 
Denote by $\kappa_{\ell}: \mathcal{A}_{1}\to\mathcal{A}_{\ell}$ the homothety of factor $2^{\ell}$, and write:
 \begin{align*}
  \int_{\R^3\backslash B} |u_0|^4\vol^{\R^3} &= \sum_{\ell=0}^{\infty} \int_{\mathcal{A}_{\ell}} |u_0|^4\vol^{\R^3} \sim \sum_{\ell=0}^{\infty}(2^{\ell})^2 \int_{\mathcal{A}_1} |\kappa_{\ell}^*u_0|^4\vol^{\R^3} \\
                                             & \leq\sum_{\ell=0}^{\infty} (2^{\ell})^2c\Big[\Big(\int_{\mathcal{A}_{1}} |d(\kappa_{\ell}^*u_0)|^2\vol^{\R^3}\Big)^2
                                                                     + \Big(\int_{\mathcal{A}_{1}} |\kappa_{\ell}^*u_0|^2\vol^{\R^3}\Big)^2\Big]\\
                                             &\sim c\sum_{\ell=0}^{\infty} \Big[\Big(\int_{\mathcal{A}_{\ell}} |du_0|^2\vol^{\R^3}\Big)^2
                                                                                +\Big(\int_{\mathcal{A}_{\ell}} |u_0|^2\rho^{-2}\vol^{\R^3}\Big)^2\Big]\\
                                             &= c\Big[\Big(\int_{\R^3\backslash B} |u_0|^2\vol^{\R^3}\Big)^2+\Big(\int_{\R^3\backslash B} |du_0|^2 \rho^{-2}\vol^{\R^3}\Big)^2\Big].
 \end{align*}
The norms here are taken for $g_{\R^3}$. 
We used the Sobolev embedding $L^{1,2}(\mathcal{A}_{1})\to L^4(\mathcal{A}_{1})$ between the first and the second lines, and denoted its norm (to the power 4) by $c$. 
\cqfd

 \subsubsection{Unweighted second order and third order estimates.}

The technique we use here is really the same as that used by Joyce in the ALE case, and which is essentially the same as in the compact case. 
It is based on the following observation.
 \begin{prop}  \label{prop_aubin_yau}
  Let $f\in C^0_{\beta+2}\cap C^2$, $\beta\in(0,1)$, and $\varphi\in C^3_{\gamma}\cap C^4_{loc}$, $0<\gamma\leq\beta$, such that $\big(\omega_Y+i\ddbar\varphi\big)^2=e^{f}\omega^2_Y$. 
Denote by $F$ the function $\log(4-\Delta\varphi)-\kappa\varphi$ on $Y$ where $\kappa$ is any constant and $\Delta$ the Laplacian of $g_Y$, and by $\Delta'$ the Laplacian operator with respect to $\omega_{\varphi}$. 
Then there exists a constant $C$ depending only on $\|\riem^{\omega_Y}\|_{C^0}$ such that
  \begin{equation*}
   \Delta'F\leq (4-\Delta\varphi)^{-1}\|\Delta f\|_{C^0}+\kappa(2-\tr^{\omega_{\varphi}}\omega)+C\tr^{\omega_{\varphi}}\omega_Y.
  \end{equation*}
 \end{prop}

We do not prove this proposition, because (up to some minor changes, like $4$, which is the real dimension of $Y$, instead of its complex dimension because we use $i\ddbar$ instead of $dd^c$) it all comes from a local formula, which is proved in \cite{yau}, \cite{aub} or \cite{joy}. 
Nevertheless, because the computation of this formula can be considered as a tour-de-force, we quote it now: in the conditions of the proposition, if $g'$ is the metric $\omega_{\varphi}(\cdot,J_Y\cdot)$, in local holomorphic coordinates and with Einstein's summation convention, one has
 \begin{equation}  \label{yau_formula}
  \begin{aligned} 
    \Delta'(\Delta \varphi)
             =&-2\Delta f+4g^{\alpha\bar{\lambda}}g'^{\mu\bar{\beta}}g'^{\gamma\bar{\nu}}\nabla_{\alpha\bar{\beta}\gamma}\varphi
               \nabla_{\bar{\lambda}\mu\bar{\nu}}\varphi\\
              &+4g'^{\alpha\bar{\beta}} g^{\gamma\bar{\delta}}\big({(\riem^{\omega_Y})^{\bar{\vareps}}}_{\bar{\delta}\gamma\bar{\beta}}
              \nabla_{\alpha\bar{\vareps}}\varphi                                                              
              -{(\riem^{\omega_Y})^{\bar{\vareps}}}_{\bar{\beta}\alpha\bar{\delta}}\nabla_{\gamma\bar{\vareps}}\varphi\big).
    \end{aligned}
 \end{equation}
Here and in what follows, $\nabla$ is the Levi-Civita connection of $g_Y$. 
We bring also the precision that $4-\Delta\varphi$ is always positive: take the trace with respect to $\omega_Y$ of $\omega_{\varphi}$, which is automatically positive as we saw when proving the uniqueness of the solution of $(E_t)$), and even $\geq 4e^{f/2}$, so that its inverse is a priori bounded above. 

Furthermore, it is not hard to deduce from Proposition \ref{prop_aubin_yau} a second order estimate on $\varphi$:
 \begin{crl}  \label{crl_scd_ord_estimate}
  Under the assumptions of Proposition \ref{prop_aubin_yau}, there exists a constant $Q\geq 0$ depending only on $\beta$, $\gamma$, $\|f\|_{C^0_{\beta+2}}$, $\|\Delta f\|_{C^0}$, $g$, $m$ and $\|\riem^{\omega_Y}\|_{C^0}$ such that $\Delta\varphi\geq -Q$.
 \end{crl}
\prf. 
Fix $\kappa=C+1$ in Proposition \ref{prop_aubin_yau}. 
Two situations can occur: the function $F$ of this proposition achieves its supremum or not. 
If not, since decay conditions imply that it tends to $2\log 2$ at infinity, we get that $4-\Delta\varphi\leq 4e^{\kappa\varphi}$, and conclude with Proposition \ref{prop_c0_estimate}.

If now $F$ achieves its supremum, at a point $p$ say, we have that $\Delta'F(p)\geq 0$, and a little computation using this inequality shows that at $p$, $\tr^{\omega'}\omega_Y\leq C':=4\kappa+\tfrac{1}{4}e^{\|f\|_{C^0}/2}\|\Delta f\|_{C^0}$. 
From this we get $F(p)\leq \|f\|_{C^0}+\log (2C')+\kappa\|\varphi\|_{C_0}$, and since $F\leq F(p)$ on $Y$, the conclusion follows, noticing that $C'$ and $\|\varphi\|_{C_0}$ only depend on the parameters announced in Propositions \ref{prop_c0_estimate} and \ref{prop_aubin_yau}. 
Details can be found in \cite{joy}. 
\cqfd

~

Using now the positivity condition on $\omega_{\varphi}$ and the upper and lower estimates on $\Delta\varphi$, one gets :
 \begin{prop}  \label{prop_scd_ord_estimate}
  Let $f\in C^0_{\beta+2}\cap C^2$, $\beta\in(0,1)$, and $\varphi\in C^3_{\gamma}\cap C^4_{loc}$, $0<\gamma\leq\beta$, such that $\big(\omega_Y+i\ddbar\varphi\big)^2=e^{f}\omega_Y^2$. 
  Then there exist some constants $Q_1$ and $Q_2$ depending only on $\beta$, $\gamma$, $\|f\|_{C^0_{\beta+2}}$, $\|\Delta f\|_{C^0}$, $g_Y$, and $\|\riem^{\omega_Y}\|_{C^0}$ such that
  \begin{equation*}
   \big\|i\ddbar\varphi\big\|_{C^0}\leq Q_1 \quad \text{and} \quad Q_2^{-1}\omega_Y\leq \omega_{\varphi} \leq Q_2\omega_Y.
  \end{equation*}
 \end{prop}

~

Finally, we give a third order estimate :
 \begin{prop}  \label{prop_thd_ord_estimate}
  Let $f\in C^0_{\beta+2}\cap C^3$, $\beta\in(0,1)$, and $\varphi\in C^3_{\gamma}\cap C^5_{loc}$, $0<\gamma\leq\beta$, such that $\big(\omega_Y+i\ddbar\varphi\big)^2=e^{f}\omega_Y^2$. 
Then there exists a constant $Q_3$ depending only on $\beta$, $\gamma$, $\|f\|_{C^0_{\beta+2}}$, $\|f\|_{C^3}$, $g$, and $\|\riem^{\omega_Y}\|_{C^1}$ such that 
$\big\|\nabla i\ddbar\varphi\big\|_{C^0}\leq Q_3$.
 \end{prop}
\prf. Here again we do not give the complete proof because it is very similar to the one in \cite{joy}. 
We only say a few words about the main ingredients. 
Set $S$ so that $4S^2=\big\|\nabla i\ddbar\varphi\big\|_{C^0}^2$, with $\nabla$ the Levi-Civita connection of $g_Y$. 
Formula \eqref{yau_formula} tells us that $\Delta'(\Delta\varphi)\geq cS^2-C$ for some constants $c>0$ and $C$ depending only on the parameters announced. 
Moreover, a hard but local computation shows that $\Delta'(S^2)$ is equal to a nonpositive quantity plus a linear term (with coefficients that are polynomials in $e^f$, $\omega_{\varphi}$, $\nabla^2f$, $\riem^{\omega_Y}$) in $\nabla i\ddbar\varphi$ and a quadratic term (with coefficients that are polynomials in $e^f$, $\omega_{\varphi}$, $\nabla^3f$, $\nabla R$) in $\nabla i\ddbar\varphi$. 
This we can sum up by saying there exists a constant $C'$ depending only on the parameters announced such that $\Delta'(S^2)\leq C'(S^2+S)$. 
Now those considerations give $\Delta'(S^2-2c^{-1}CC'\Delta\varphi)\leq -C\big(S-\tfrac{1}{2}\big)^2+2c^{-1}CC'+\tfrac{1}{4}C$; 
one concludes according to $(S^2-2c^{-1}CC'\Delta\varphi)$ achieves its supremum at some point (and using its $\Delta'$ is $\geq0$ at this point) or not (and using decay conditions).
\cqfd

 \subsubsection{$C^0_{\delta}$ estimates.}  \label{c0delta_estimate}

We now prove an estimate on a weighted $C^0$-norm of a solution of a complex Monge-Ampère equation on $Y$. 

For this, we take a point of view close to that of paragraph \ref{c0_estimate} but we estimate $\|\varphi\rho^{\delta}\|_{C^0}$ (i.e. $\|\varphi\|_{C^0_{\delta}}$), $\delta\in(0,\gamma)$, instead of $\|\varphi\|_{C^0}$ with the same iteration. 
Of course, putting a weight $\rho^{p\delta}$ makes the integrations by parts a bit more complicated, but the whole spirit of the proof remains the same. 
Here again we write the computations in detail since they differ at some points from Joyce's. 

To begin with, we state the weighted integration by parts formula we are using in our iteration scheme:
 \begin{lem}  \label{lem_weighted_ibp_formula}
  Let $f\in C^0_{\beta+2}$, $\beta\in(0,1)$, and $\varphi\in C^3_{\gamma}$, $0<\gamma\leq\beta$, such that $\big(\omega_Y+i\ddbar\varphi\big)^2=e^{f}\omega_Y^2$. 
  Set $T=\omega_Y+\omega_{\varphi}$. 
  Then for any $p>2$, $q$ such that $p\gamma-q>1$, 
  \begin{align*}
   \int_Y \big|\partial\big(|\varphi|^{p/2}\rho^{q/2}\big)&\big|^2\omega^2_Y \leq \frac{p^2}{2(p-1)}\int_Y \rho^q|\varphi|^{p-2}\varphi (e^f-1)\omega^m_Y\\
                                                          & + \frac{q}{2(p-1)}\int_Y |\varphi|^p\rho^{q-2}\big[(p+q-2)i\partial\rho\wedge\dbar\rho-(p-2)\rho i\ddbar\rho\big]\wedge T,
  \end{align*}
 where those integrals are finite. 
 \end{lem}
\prf. As in \cite{joy}, we say that $\int_Y d\big[p^2\varphi|\varphi|^{p-2}\rho^qd^c\varphi\wedge T+(p-2)q|\varphi|^p\rho^{q-1}\wedge T\big]=0$, because of the same Stokes' argument as in the proof of Lemma \ref{lem_ibp_formula}. 
Now, a direct computation yields 
 \begin{align*}
  d\big[p^2\varphi&|\varphi|^{p-2}\rho^qd^c\varphi\wedge T+(p-2)q|\varphi|^p\rho^{q-1}\wedge T\big]\\
       = & p^2(p-1)|\varphi|^{p-2}\rho^qd\varphi\wedge d^c\varphi\wedge T+p^2|\varphi|^{p-2}\rho^q dd^c\varphi\wedge T \\
         & +2pq(p-1)\varphi|\varphi|^{p-2}\rho^{q-1}d^c\varphi\wedge T+q(q-1)(p-2)|\varphi|^{p}\rho^{q-2}d\rho\wedge d^c\rho\wedge T\\
         & +(p-2)q|\varphi|^p \rho^{q-1}dd^c\rho\wedge T
 \end{align*}
(notice that $d\varphi\wedge d^c\rho\wedge T=d\rho\wedge d^c\varphi\wedge T$ since $T$ is of bidegree $(1,1)$). 
On the other hand,
 \begin{align*}
  d\big(|\varphi|^{p/2}\rho^{q/2}\big)\wedge d^c\big(|\varphi|^{p/2}\rho^{q/2}\big)\wedge T 
       = &  \frac{p^2}{4}|\varphi|^{p-2}\rho^q d\varphi\wedge d^c\varphi +\frac{q^2}{4}|\varphi|^p\rho^{q-2}d\rho\wedge d^c\rho\\
         & + \frac{pq}{2} |\varphi|^{p-2}\varphi\rho^{q-1}d\rho\wedge d^c\varphi\wedge T.
 \end{align*}
Comparing those identities one can write
 \begin{align*}
  d\big[p^2\varphi|\varphi|^{p-2}&\rho^qd^c\varphi\wedge T+(p-2)q|\varphi|^p\rho^{q-1}\wedge T\big]\\
       = & 4(p-1) d\big(|\varphi|^{p/2}\rho^{q/2}\big)\wedge d^c\big(|\varphi|^{p/2}\rho^{q/2}\big)\wedge T+p^2\varphi|\varphi|^{p-2}\rho^q dd^c\varphi\wedge T \\
         & +q|\varphi|^p\rho^{q-2}\big[(p-2)\rho dd^c\rho\wedge T-(p+q-2)d\rho\wedge d^c\rho\wedge T\big]. 
 \end{align*}
Then conclude after dividing this by 2, integrating both sides over $Y$ and noticing that $i\ddbar\varphi\wedge T= (e^f-1)\omega^2$ and $T\geq \omega_Y$. 
Checking that all resulting integrals are finite is straightforward with the assumption $p\gamma-q>1$. 
\cqfd

~

Now fix $\delta\in(0,\gamma)$; we take $q$ in the latter lemma as $p\delta$, so the condition $p\gamma-q>1$ becomes $p>\tfrac{1}{\gamma-\delta}$. 
The following lemma initiates our sequence of recursive controls: 
 \begin{lem}  \label{lem_initiate_weighted}
  Fix $p_0'>2$, $p_0'>\tfrac{1}{\gamma-\delta}$, $p'_0\geq p_0=\tfrac{2}{\gamma}$. 
  Under the assumptions of Proposition \ref{prop_scd_ord_estimate}, there exists a constant $C_0$ depending only on $p'_0$, $\beta$, $\gamma$, $\delta$, $\|f\|_{C^0_{\beta+2}}$, $\|f\|_{C^2}$, $g$, and $\|\riem^{\omega_Y}\|_{C^0}$ 
such that $\|\varphi\rho^{\delta}\|_{L^{p'_0}_{d\lambda}}\leq C_0$. 
 \end{lem}
We recall that $d\lambda$ is the measure $\rho^{-1}\vol^{g_Y}$. 

~

\prf. From Lemma \ref{lem_weighted_ibp_formula} and Proposition \ref{prop_scd_ord_estimate}, and also the facts that $d\rho=O(1)$ and $dd^c\rho=O(\rho^{-1})$, we get two constants $c_1$, $c_2$ depending only on the parameters announced such that for every $p>\tfrac{1}{\gamma-\delta}$,
 \begin{equation} \label{eqn_c1c2}
  \int_Y \big|\partial\big(|\varphi|^{p/2}\rho^{q/2}\big)\big|^2\vol^{g_Y}\leq c_1p\int_Y |\varphi|^{p-1}\rho^{p\delta-(\beta+2)}\vol^{g_Y} + c_2 p\int_{Y}|\varphi|^p\rho^{p\delta-2}\vol^{g_Y}.
 \end{equation}
We now play the same game of Hölder inequalities as in paragraph \ref{c0_estimate} for the first summand of the right-hand side. 
Indeed,
 \begin{equation*}
  \int_Y |\varphi|^{p-1}\rho^{p\delta-(\beta+2)}\vol^{g_Y}\leq \Big(\int_Y \rho^{bt}\vol^{g_Y}\Big)^{1/t}\Big(\int_Y (|\varphi|\rho^{\delta})^{2p} \rho^{-1}\vol^{g_Y}\Big)^{(p-1)/2p}.
 \end{equation*}
with $t=\big(1-\tfrac{p-1}{2p}\big)^{-1}=\tfrac{2p}{p+1}$ and $b=\tfrac{p-1}{2p}+\delta-(\beta+2)$. 
Notice that the condition $bt<-3$, which ensures the convergence of $\int_Y \rho^{bt}\vol^{g_Y}$, is equivalent to $p>\tfrac{1}{\beta-\delta}$, and is thus automatically verified for the $p$ we work with. 
Furthermore, when $p$ goes to $\infty$, $bt$ tends to $-3+2(\delta-\beta)$, which is $<-3$; the quantity $\big(\int_Y \rho^{bt}\vol^{g_Y}\big)^{1/t}$ is thus bounded above by some constant $C_1$ depending only on our parameters for the considered $p$, if we choose them away from $p_0=\tfrac{2}{\gamma}$; the choice $p\geq p'_0$ is convenient.

On the other hand, we have to argue in a slightly different way for the summand $\int_Y|\varphi|^{p}\rho^{p\delta-2}\vol^{g_Y}$ of \eqref{eqn_c1c2}. 
Let us write 
 \begin{equation*}
  \int_Y|\varphi|^{p}\rho^{p\delta-2}\vol^{g_Y}  
   \leq \Big(\int_Y |\varphi|^r\rho^{-1}\vol^{g_Y}\Big)^{1/t}\Big(\int_Y (|\varphi|\rho^{\delta})^{2p} \rho^{-1}\vol^g\Big)^{1/s},
 \end{equation*}
where $\tfrac{1}{s}+\tfrac{1}{t}=1$, $\tfrac{2p\delta-1}{s}-\tfrac{1}{t}=p\delta-2$ and $\tfrac{2p}{s}+\tfrac{q}{t}=p$. 
This gives $s=\tfrac{2p\delta}{p\delta-1}:=s(p)$ (well-defined and $>2$ as soon as $p>\tfrac{1}{\delta}$, and tends to $2$ as $p$ goes to $\infty$), $t=\tfrac{2p\delta}{p\delta+1}$ (which tends to $2$) and $r=\tfrac{t}{\delta}=\tfrac{2p}{p\delta+1}$ (which tends to $\tfrac{2}{\delta}$). 
Now the condition $r\gamma>2$ ensures the convergence of $\int_Y |\varphi|^r\rho^{-1}\vol^{g_Y}$, and this condition turns out to be equivalent to $p>\tfrac{1}{\gamma-\delta}$, which we assume. 

Moreover, when $p$ ranges over $[p_0',\infty)$ for $p_0'>\tfrac{1}{\gamma-\delta}$ (e.g. $p_0'=\tfrac{2}{\gamma-\delta}$), $r$ ranges over $[r_0,\tfrac{2}{\delta})$, with $r_0=\tfrac{2p'_0}{p'_0\delta+1}>2$. 
From Lemma \ref{lem_induction}, $\big(\int_Y |\varphi|^r\rho^{-1}\vol^{g_Y}\big)^{1/t}$ (converges and) is bounded by some constant $C_2$. 
This constant $C_2$ depends only on the parameters announced. 
Finally, we can sum all this up saying that for $p\geq p_0'$, 
 \begin{equation*}
  \int_Y \big|\partial\big(|\varphi|^{p/2}\rho^{q/2}\big)\big|^2\vol^{g_Y}
      \leq c_1C_1p \|\varphi\rho^{\delta}\|_{L^{2p}_{d\lambda}}^{p-1}+ c_2C_2p \|\varphi\rho^{\delta}\|_{L^{2p}_{d\lambda}}^{2p/s(p)}.
 \end{equation*}
 
Take $p=p'_0$, and apply inequality \eqref{eqn_sob_inj} to the LHS, with $u=|\varphi|^{p/2}\rho^{q/2}$; this yields :
 \begin{equation*}
  \|\varphi\rho^{\delta}\|_{L^{2p_0'}_{d\lambda}}^{p_0'}  
      \leq C_{S}c_1C_1p_0' \|\varphi\rho^{\delta}\|_{L^{2p_0'}_{d\lambda}}^{p_0'-1}
          + C_{S}c_2C_2p_0' \|\varphi\rho^{\delta}\|_{L^{2p_0'}_{d\lambda}}^{2p_0'/s(p_0')},
 \end{equation*}
and one concludes noticing that $p_0'>p'_0-1$ and $p'_0>2p_0'/s(p_0')$, as $s(p_0')>2$. 
\cqfd

~

We fix $p_0'=\tfrac{2}{\gamma-\delta}$, so that it verifies all the assumptions of the latter lemma. 
As in paragraph \ref{c0_estimate}, the last step makes the transition from $\|\varphi\rho^{\delta}\|_{L^{p}_{d\lambda}}$ to $\|\varphi\rho^{\delta}\|_{L^{2p}_{d\lambda}}$ for big enough $p$ :
 \begin{lem}  \label{lem_propagate_weighted}
  Under the assumptions of Proposition \ref{prop_scd_ord_estimate}, there exists constants $C_1$ and $C_2$ only depending on $\beta$, $\gamma$, $\delta$, $\|f\|_{C^0_{\beta+2}}$, $\|f\|_{C^2}$, $g$, $m$ and $\|R^{\omega}\|_{C^0}$ 
such that for all $p\geq p'_0=\tfrac{2}{\gamma-\delta}$, $\|\varphi\rho^{\delta}\|_{L^{2p}_{d\lambda}}^p\leq C_1p\|\varphi\rho^{\delta}\|_{L^{p}_{d\lambda}}^{p-1} + C_2p\|\varphi\rho^{\delta}\|_{L^{p}_{d\lambda}}^p$.
 \end{lem}
\prf. We saw in the beginning of the proof of Lemma \ref{lem_initiate_weighted} that there exist $c_1$ and $c_2$ depending only on the parameters announced such that for all $p>\tfrac{1}{\gamma-\delta}$
 \begin{equation*}
  \int_Y \big|\partial\big(|\varphi|^{p/2}\rho^{q/2}\big)\big|^2\vol^{g_Y}  \leq c_1p\int_Y |\varphi|^{p-1}\rho^{p\delta-(\beta+2)}\vol^{g_Y} 
                                                                           + c_2 p\int_Y|\varphi|^p\rho^{p\delta-2}\vol^{g_Y}.
 \end{equation*}
We deal again with the two summands of the right-hand side separately. 
For the concerned $p$, $\int_Y|\varphi|^{p-1}\rho^{p\delta-(\beta+2)}\vol^{g_Y}=\int_Y(|\varphi|\rho^{\delta})^{p-1}\rho^{\delta-(\beta+2)}\vol^{g_Y}$ and by Hölder inequality this is bounded above by 
 \begin{equation*}
  \Big(\int_Y\rho^{p(\delta-\beta)-p-1}\vol^{g_Y}\Big)^{1/p} \Big(\int_Y(|\varphi|\rho^{\delta})^p\rho^{-1}\vol^{g_Y}\Big)^{(p-1)/p}.
 \end{equation*}
The first integral converges for $p>\tfrac{2}{1+\beta-\delta}$, which is automatic if $p>\tfrac{1}{\gamma-\delta}$, 
and its $\tfrac{1}{p}$-th power tends to $\sup_Y \rho^{\delta-\beta-1}$ as $p$ goes to infinity, so we get the constant $C_1$ of the statement (after multiplying by $c_1$), if we restrict to $p\in[\tfrac{2}{\gamma-\delta},\infty)$. 

Let us deal next with the summand $\int_Y|\varphi|^p\rho^{p\delta-2}\vol^{g_Y}$; we can take $C_2=c_2\sup_{Y}\rho^{-1}$ without more efforts, and the lemma is proved. 
\cqfd

~

Under the assumptions of Proposition \ref{prop_scd_ord_estimate} and the condition $\delta\in(0,\gamma)$, it is again an easy exercise, in view of Lemmas \ref{lem_initiate_weighted} and \ref{lem_propagate_weighted}, 
to get two constants $Q_{\delta}$ and $C_3$ depending only on the same parameters as those announced in Lemma \ref{lem_propagate_weighted} such that for every $p\geq p_0'$, we have $\|\varphi\rho^{\delta}\|_{L^{p}_{d\lambda}}\leq Q_{\delta}(C_3p)^{-2/p}$. 
Letting $p$ go to $\infty$, we can conclude:
 \begin{prop}  \label{prop_c0_delta_estimate}
  Under the assumptions of Proposition \ref{prop_scd_ord_estimate} and assuming $\delta\in(0,\gamma)$, there exists a constant $Q_{\delta}$ depending only on $\beta$, $\gamma$, $\delta$, $\|f\|_{C^0_{\beta+2}}$, $\|f\|_{C^2}$, $g_Y$, and $\|\riem^{\omega_Y}\|_{C^0}$ such that $\|\varphi\|_{C^0_{\delta}}\leq Q_{\delta}$.
 \end{prop}

 \subsubsection{$C^{k,\alpha}$ and $C^{k,\alpha}_{\delta}$ estimates.}

We shall now state (unweighted) higher order a priori estimates :

 \begin{prop}  \label{prop_c_k_alpha_estimate}
  Let $f\in C^0_{\beta+2}\cap C^{k,\alpha}$, $\beta\in(0,1)$, $k\geq 3$, $\alpha\in(0,1)$ and $\varphi\in C^3_{\gamma}\cap C^5_{loc}$, $0<\gamma\leq\beta$, such that $\big(\omega_Y+i\ddbar\varphi\big)^2=e^{f}\omega_Y^2$. 
Then $\varphi\in C^{k+2,\alpha}$ and there exists a constant $Q_{k,\alpha}$ only depending on $\beta$, $\gamma$, $k$, $\alpha$, $\|f\|_{C^0_{\beta+2}}$, $\|f\|_{C^{k,\alpha}}$ and $\|g_y\|_{C^{k,\alpha}}$ such that 
$\|\varphi\|_{C^{k+2,\alpha}}\leq Q_{k,\alpha}$.
 \end{prop}
\prf. The theorem follows from an inductive use of the crucial Aubin-Yau formula \eqref{yau_formula}, which we recall here: 
 \begin{align*}
  \Delta'(\Delta \varphi)
             =&-2\Delta f+4g^{\alpha\bar{\lambda}}_Yg'^{\mu\bar{\beta}}g'^{\gamma\bar{\nu}}\nabla_{\alpha\bar{\beta}\gamma}\varphi
               \nabla_{\bar{\lambda}\mu\bar{\nu}}\varphi\\
              &+4g'^{\alpha\bar{\beta}} g^{\gamma\bar{\delta}}\big({(\riem^{\omega_Y})^{\bar{\vareps}}}_{\bar{\delta}\gamma\bar{\beta}}
              \nabla_{\alpha\bar{\vareps}}\varphi                                                              
              -{(\riem^{\omega_Y})^{\bar{\vareps}}}_{\bar{\beta}\alpha\bar{\delta}}\nabla_{\gamma\bar{\vareps}}\varphi\big)
 \end{align*}
($\Delta'$ stands for the scalar Laplacian with respect to $\omega_{\varphi}$).

To begin with, suppose $k=3$. 
Then according to the assumptions on $f$ and the conclusions of Propositions \ref{prop_scd_ord_estimate} and \ref{prop_thd_ord_estimate}, we have for the moment a $C^0$ estimate in terms of the parameters announced on the right-hand side. 

Now the geometry of $Y$ and the asymptotics of $g_Y$ allow us to take an atlas of $Y$ of holomorphic balls with uniform radius, such that the pull-backs of $g_Y$ to these balls are uniformly bounded above and below, and so are the pull-backs of $g'$. 
As a consequence the Laplacian operators associated to these latter Kähler metrics are uniformly elliptic in $C^0$ sense. 
We can ask furthermore that the family of those balls with the half radius still gives an atlas.

To say this more precisely, we have a family of holomorphic charts $\pi_i: B(0,1)\rightarrow Y$, $i\in \mathcal{I}$, such that $Y=\bigcup_{i\in \mathcal{I}}\pi_i\big(B(0,1)\big)=\bigcup_{i\in \mathcal{I}}\pi_i\big(B(0,\tfrac{1}{2})\big)$ and such that there exist constants $c>0$ and $C$ depending on the parameters announced such that for all $i\in \mathcal{I}$, $c \e\leq \pi_i^*g'\leq C \e$. 
Then the $\pi_i^*\Delta'$ are uniformly elliptic, in the $C^0$ sense.

Pulling back the formula above by any $\pi_i$, the right-hand side is bounded independently of $i$, and so is $\pi_i^*(\Delta \varphi)$. 
The standard Morrey-Schauder's theorem tells us that the $\pi_i^*(\Delta \varphi)$ is bounded in $C^{1,\alpha}$ on $B(0,\tfrac{3}{4})$, again independently of $i$; a careful reading shows quickly that this uniform bound depends only on the parameters of the statement. We can reformulate all this saying that $\Delta\varphi$ is in $C^{1,\alpha}$, and the corresponding norm is controlled in terms of the parameters.

On the other hand, $\pi_i^*(\Delta\varphi)=(\pi_i^*\Delta)(\pi_i^*\varphi)$. 
We can also suppose our covering is taken so that there exist constants $c>0$ and $C$ depending only on the parameters announced such that for all $i\in \mathcal{I}$, $c \e \leq \pi_i^*g$ and $\|g\|_{C^{1,\alpha}}\leq C$. 
This allows us to apply again a Schauder estimate, and conclude that we have a uniform bound on the $\pi_i^*\varphi$ in $C^{3,\alpha}\big(B(0,\tfrac{1}{2})\big)$, that is: $\varphi\in C^{3,\alpha}$ and the corresponding norm is again controlled in terms of the parameters announced.

This settles the $k=3$ case. 
Now observe that a $C^{3,\alpha}$ bound on $\varphi$ together with the assumptions on $f$  (at least, a $C^{2,\alpha}$ bound) give a $C^{0,\alpha}$ bound on the right-hand side of \eqref{yau_formula}. 
Applying again twice the Schauder estimates (once with $\Delta'$ which is uniformly elliptic in the $C^{0,\alpha}$ sense from the previous case, and once to $\Delta$ which is uniformly elliptic in the $C^{2,\alpha}$ sense) after refining the covering if needed so, one gets the announced $C^{4,\alpha}$ estimate on $\varphi$. 
Following once more this scheme using in particular the $C^{k,\alpha}$ bounds on $f$, and noticing the Schauder estimates include regularity before giving bounds on $C^{\ell,\alpha}$ norms, we get the result with $k=5$. 
For the general case, just follow the initiated induction. 
\cqfd

~

We state a similar result for the \textit{weighted} higher order estimates :

 \begin{prop}  \label{prop_c_k_alpha_gamma_estimate}
  Let $f\in C^{k,\alpha}_{\beta+2}\cap C$, $\beta\in(0,1)$, $k\geq 3$, $\alpha\in(0,1)$ and $\varphi\in C^3_{\gamma}\cap C^5_{loc}$, $0<\gamma\leq\beta$, such that $\big(\omega_Y+i\ddbar\varphi\big)^2=e^{f}\omega_Y^2$. 
Then $\varphi\in C^{k+2,\alpha}_{\gamma}$. 
Moreover, there exists a constant $Q_{k,\alpha,\gamma}$ depending on only $\beta$, $\gamma$, $k$, $\alpha$, $\|\varphi\|_{C^{k,\alpha}_{\gamma}}$ and $\|f\|_{C^{k,\alpha}_{\beta}}$ 
such that 
$\|\varphi\|_{C^{k+2,\alpha}_{\gamma}}\leq Q_{k,\alpha,\gamma}$.
 \end{prop}
\prf. We use a technical lemma from \cite{joy}, (Proposition 8.6.12, p.197), which remains true in our ALF setting; its utility lies in the rescaling process used in establishing elliptic weighted estimates. 
The lemma states, with our notations:
 \begin{lem}  \label{joyce_technical_prop}
  Let $K_1$, $K_2>0$, $\lambda\in[0,1]$ and $k\geq 3$ be an integer. 
  Then there exists $K_3$ depending only on $\alpha$, $\beta$, $\gamma$, $\|\varphi\|_{C^{0}_{\gamma}}$ and $K_1,K_2,\lambda,k$ such that the following holds:
  
  Under the assumptions of Proposition \ref{prop_c_k_alpha_gamma_estimate} and if $\|f\|_{C^{k,\alpha}_{\beta}}\leq K_1$, 
  $\big\|\nabla^jdd^c\varphi\big\|_{C^0_{\lambda\ell}}\leq K_2$, $\ell=0,\dots,k$ and $\big[\nabla^kdd^c\varphi\big]^{\alpha}_{\lambda k+(\lambda-1)\alpha}\leq K_2$, then:
  $\big\|\nabla^jdd^c\varphi\big\|_{C^0_{\lambda\ell+\gamma}}\leq K_3$, $\ell=0,\dots,k+2$ and $\big[\nabla^kdd^c\varphi\big]^{\alpha}_{\gamma+\lambda(k+2)+(\lambda-1)\alpha}\leq K_3$.
 \end{lem}
The $[\cdot]^{\alpha}_{\cdots}$ are analogous to those used in \ref{df_wtd_hldr_sp}, with $Y$, $g_Y$ and $\rho$ instead of $X$, $g_m$ and $\tilde{R}$. 

~

\prf \textit{of Lemma \ref{joyce_technical_prop}}. 
The major change here is that the injectivity radius does not grow as fast as $\rho$, but instead remains bounded, essentially by half the length of the fibers. 
But this is not an issue. 
Indeed, the Riemannian exponential map is still well defined, and authorizes the following manipulations. 
Given $x\in \{\rho\geq 2\rho_0\}$ ($\rho_0$ determined later), identify $(T_xY,J_Y,g_Y)$ with and $(\C^2,I_1,\e)$. 
Take $R>0$, and denote by $\pi_{R,x}$ the map $B_{\e}(0,1)\rightarrow B_{g_Y}(x,R)$, $y\mapsto \exp_x^{g_Y}(Ry)$; it is not a diffeomorphism in general, since large balls wrap following asymptotically the fibers of $\varpi$. 
We can nonetheless define the operator $P_{x,R}:C^{k+2,\alpha}\big(B_{\e}(0,1)\big)\rightarrow C^{k,\alpha}\big(B_{\e}(0,1)\big)$ by
 \begin{equation*}
  P_{x,R}(v)=R^2 \frac{\big({\pi_{R,x}}^*(i\ddbar)\big)(v)\wedge{\pi_{R,x}}^*(\omega_Y+\omega_{\varphi}^{m-1})}{{\pi_{R,x}}^*(\omega_Y^2)}.
 \end{equation*}
Then one can take $R=L\rho(x)^{\lambda}$ with $L=L(\rho_0,\lambda,g_Y)$ small enough so that $B_{g_Y}(x,R)\subset\{\rho\geq\rho_0\}$; this way one has
 \begin{align*}
  &\big\|R^{-2}{\pi_{x,R}}^*g_Y - \e \big\|_{C^{k,\alpha}}                \leq \frac{1}{2} \quad\text{for all }x\in\{\rho\geq2\rho_0\},\quad\text{and} \\ 
  &\big\|R^{-2}{\pi_{x,R}}^*\omega_Y - \omega_{\e}\big\|_{C^{k,\alpha}}  \leq \frac{1}{2} \quad\text{for all }x\in\{\rho\geq2\rho_0\},  
 \end{align*}
if $\rho_0$ is chosen big enough, thanks to the asymptotic geometry of $g_Y$. 
Now the rest of Joyce's proof applies unchanged (in particular, one is brought to using Schauder estimates between the \textit{fixed} balls $B_{\e}(0,2)$ and $B_{\e}(0,1)$, with a $C^{k,\alpha}$ uniformly elliptic family of operators), since the identity 
 \begin{equation*}
  P_{x,R}\big((\pi_{x,R})^*\varphi\big)= R^2\big(e^{(\pi_{x,R})^*f}-1\big)
 \end{equation*}
is again just a rewriting of the pulled-back Monge-Ampère equation verified by $\varphi$. 

To complete the present proof though, one has to deal with the compact subset $\{\rho\leq 2\rho_0\}$; 
in that case, the estimates of the statement are an immediate consequence of the unweighted estimates of Proposition \ref{prop_c_k_alpha_estimate}. 
\hfill $\blacksquare$

~

Now Proposition \ref{prop_c_k_alpha_gamma_estimate} follows from a repeated use of Lemma \ref{joyce_technical_prop}, as in \cite[p. 195]{joy}. 
At the end of the argument, we pass to $\lambda=1$ in this lemma, which precisely gives estimates on the $\big\|\nabla^jdd^c\varphi\big\|_{C^0_{\ell+\gamma}}$, $\ell=0,\dots,k+2$, 
and on $\big[\nabla^kdd^c\varphi\big]^{\alpha}_{\gamma+(k+2)+\alpha}$, i.e. on $\|dd^c\varphi\|_{C^{k,\alpha}_{\gamma}}$. 
\cqfd

  \subsubsection{Closedness of $S$}   \label{prgph_clsdnss}

  We conclude this section by the proof of the closedness of $S$. 
Take a sequence $(t_j)$ of elements of $S$, converging to some $t_{\infty}\in [0,1]$; in particular, each $\varphi_{t_j}\in C^{\infty}_{\beta/2}$.  
Pick $\alpha\in(0,1)$. 
From Proposition \ref{prop_c_k_alpha_gamma_estimate}, since the $t_jf$ now play the role of the $f$ of that Proposition, we have a bound on $\|\varphi_{t_j}\|_{C^{5,\alpha}_{\beta/2}}$  which depends only on $\|f\|_{C^{3,\alpha}_{\beta}}$ and the $\|\varphi_{t_j}\|_{C^{0}_{\beta/2}}$. 
But from Proposition \ref{prop_c0_delta_estimate}, we also have a bound on the $\|\varphi_{t_j}\|_{C^{0}_{\beta/2}}$, which depends only on $\|f\|_{C^{3,\alpha}_{\beta}}$. 
In other words, $(\varphi_{t_j})$ is bounded in $C^{5,\alpha}_{\beta/2}$, and thus converges weakly to some $\varphi_{\infty}\in C^{5,\alpha}_{\beta/2}$. 
Now the inclusion $C^{5,\alpha}_{\beta/2}\hookrightarrow C^{5,\alpha/2}_{\beta/4}$ is compact; we can thus assume that $(\varphi_{t_j})$ converges strongly to $\varphi_{\infty}$ in $C^{5,\alpha/2}_{\beta/4}$. 
In particular the $(E_{t_j})$ pass to the limit to give:
 \begin{equation*}
  \big(\omega_{Y}+i\ddbar\varphi_{\infty}\big)^2=e^{t_{\infty}f}\omega_Y^2.
 \end{equation*}
Developing this latter equation yields: $\omega_Y^2 +2\omega_Y\wedge i\ddbar\varphi_{\infty}+(i\ddbar\varphi_{\infty})^2=e^{f_{t_{\infty}}}\omega_Y^2$, 
that is: 
 \begin{equation*}
  -\frac{1}{2}\Delta\varphi_{\infty}= \big(e^{t_{\infty}f}-1\big)+\frac{(i\ddbar\varphi_{\infty})^2}{\omega_Y^2} \in C^{3,\alpha/2}_{\beta+2}.
 \end{equation*}
On the other hand, $\Delta: C^{5,\alpha/2}_{\gamma}\to C^{3,\alpha/2}_{\gamma+2}$ is an isomorphism for any $\gamma\in (0,1)$; this follows from \cite[App. A]{bm}. 
In view of the latter equation, we hence get $\varphi_{\infty}\in C^{5,\alpha/2}_{\beta}$. 
Apply Proposition \ref{prop_c_k_alpha_gamma_estimate} to conclude that $\varphi_{\infty}\in C^{\infty}_{\beta}$, and thus that $t_{\infty}\in S$. 
\cqfd

\subsection{Openness of $S$.}   \label{sect_openness}

The only point we are left with is the openness of $S$; we settle it now. 
Let $t\in S$; write $\omega_t$ for $\omega_Y+i\ddbar\varphi_t$. 
We already know that $\omega_t$ is Kähler. 
Moreover, the linearisation of the Monge-Ampère operator
 \begin{align*}
  C^{5,\alpha}_{\beta} &\longrightarrow C^{3,\alpha}_{\beta+2} \\
        \psi          &\longmapsto     \frac{\big(\omega_Y+i\ddbar\psi\big)^2}{(\omega_Y)^2},
 \end{align*}
with $\alpha\in (0,1)$, at point $\varphi_t$, is $-\tfrac{e^{tf}}{2}\Delta_t$, where $\Delta_t$ is the Laplacian of $\omega_t$. 
Now $\omega_t$ is $C^{\infty}_{\beta+2}$ close to $\omega_Y$; it follows that $\Delta_t$, and hence $-\tfrac{e^{tf}}{2}\Delta_t$, are isomorphisms $C^{5,\alpha}_{\beta} \to C^{3,\alpha}_{\beta+2}$. 
Consequently, $(E_s)$ has a (necessarily unique) solution $\varphi_s\in C^{5,\alpha}_{\beta}$ for all $s$ close to $t$. 
At last, we apply Proposition \ref{prop_c_k_alpha_gamma_estimate} to get that those $\varphi_s$ are in $C^{\infty}_{\beta}$, 
which ends the proof of the openness of $S$, and that of Theorem \ref{thm_CY_ALF}. 
\cqfd

~


\begin{small}

 \renewcommand{\refname}{\large{References}}

\end{small}

~

\small \textsc{École Normale Supérieure, UMR 8553}

\url{hugues.auvray@ens.fr}

\end{document}